\documentclass[10pt,a4paper]{article}
\usepackage[utf8]{inputenc}
\usepackage{algorithm}
\usepackage{algorithmic}
\usepackage {amsmath}
\usepackage{amssymb}
\usepackage{amsfonts}
\usepackage{amsthm}
\usepackage{ascmac}
\usepackage{color}
\usepackage{framed}
\newcommand\norm[1]{\left\lVert#1\right\rVert}
\newtheorem{df.}{Definition}
\newtheorem{th.}{Theorem}
\newtheorem{lm.}{Lemma}
\newtheorem{cd.}{Condition}
\newtheorem{ex.}{Example}
\begin {document}
\begin{titlepage}
\title { Bayesian Generalization Error of  Poisson Mixture and Simplex Vandermonde Matrix Type Singularity}
\author{Kenichiro SATO\footnote{E-mail: sato.k.bo@m.titech.ac.jp}  
{ and }Sumio WATANABE\footnote{E-mail: swatanab@c.titech.ac.jp}\\
Department of Mathematical and Computing Science\\
 Tokyo Institute of Technology\\
2-12-1, Oookayama, Meguro-ku, Tokyo, 152-8552, Japan}
\date{}
\end{titlepage}

\maketitle
\begin{abstract}
A Poisson mixture is one of the practically important models in computer science, biology, and sociology. However, the theoretical property has not been studied because the posterior distribution can not be approximated by any normal distribution. Such a model is called singular and it is known that Real Log Canonical Threshold (RLCT) is equal to the coefficient of the asymptotically main
term of the Bayesian generalization error. 
In this paper, we derive RLCT of a simplex Vandermonde matrix type singularity which
 is equal to that of a Poisson mixture in general cases.
\end{abstract}

\section{Introduction}

In this section, Bayesian inference is introduced and several notations are defined, 
which are summarized in Table.\ref{table:111}. 

Let $q(x)$ be a probability density function defined on $\mathbb{R}^M$ and 
$X^n$ be a set of random variables which are independently subject to $q(x)$, 
where  $X^n$ and $n$ are called a sample and a sample size, respectively. 
A statistical model $p(x|w)$ and a prior are defined by a conditional 
probability density function of $x$ for a given parameter $w$ and 
a probability density function $\varphi (w)$, respectively, where 
$w\in W\subset\mathbb{R}^d$. In this paper, we assume that the set of parameters 
$W$ is a sufficiently large compact set whose interior is not the empty set. 
The posterior distribution is defined by
\begin{eqnarray*}
p(w|X^n) &:=& \frac {1} {Z_n } \varphi (w) \prod_{i=1}^n p(X_i | w ), 
\end{eqnarray*}
where $Z_n := \int_W \varphi(w) \prod_{i=1}^n p(X_i | w )  dw $ is a normalizing constant. 
The free energy $F_n$ is also defined by
$F_n= - {\rm log} Z_n$.
The expected value of a function $f(w)$ by the posterior distribution $p(w|X^n)$ is 
referred to as 
\begin{eqnarray*}
E_w [ f(w)] := \int_W f(w) p(w | X^n ) dw.
\end{eqnarray*}
The Bayesian predictive distribution is defined by $E_w[p(x|w)]$ and its 
generalization error $G_n$ is given by 
\begin{eqnarray*}
G_n = - \int q(x)\log E_w [p(x|w)] dx .
\end{eqnarray*}
In this paper, we assume that there exists a parameter $w_0$ in the interior of W which minimizes Kullback-Leibler divergence 
\begin{eqnarray*}
w_0 := \arg \min_{w \in W} \left(\int q(x) \log \frac {q(x)} {p(x\mid w)} dx\right). 
\end{eqnarray*}
Note that such a parameter $w_0$ is not unique in general. We assume that $p(x|w_0)$ 
does not depend on choice of $w_0$. 
The mean error function $K(w)$ is defined by
\begin{eqnarray*}
K(w) := - \int q(x) {\rm log}  \frac {p(x|w)} {p(x|w_0)} dx.
\end{eqnarray*}

\begin{table}[htb]\label{table:111}
 \begin{center}
  \begin{tabular}{|l|c|} \hline
    Symbol & Meaning \\ \hline
	$n$ &Sample size\\
	$x$ &M-dimensional vector\\
	$w$ &parameter\\
	$W$ &set of parameters\\
	$q(x)$ & true probability distribution\\
	$p(x|w)$ & statistical model\\
	$\varphi (w)$& prior  \\
     $X^n$ & sample i.i.d. from $q(x)$\\
	$ [m: n] $ &$ \{a \in \mathbb{Z} \mid m \le a \le n \}$\\
	$ \mathbb{Z}_{\ge 0} $ &$\{a \in \mathbb{Z} \mid a \ge 0 \}$\\
	$ \mathbb{R}_{\ge 0} $ &$\{a \in \mathbb{R} \mid a \ge 0 \}$\\ 
	$ \mathbb{N} $ &$\{a \in \mathbb{Z} \mid a > 0 \}$\\
	$ \mathbb{R}_{> 0} $ &$\{a \in \mathbb{R} \mid a > 0 \}$\\
     $\delta_{i,j}$ & $1$ if $i=j$, else $0$ (Kronecker delta) \\
\hline
  \end{tabular}
\caption{Notations used in this paper}
\end{center}
\end{table}

One of the most important purposes of Bayesian statistical theory is to clarify
the asymptotic form of the Bayesian generalization error for a given triple $(q(x),p(x|w),\varphi(w))$. 
In 2000, learning theory that studies Bayesian generalization error of singular learning models was established by Watanabe \cite{watanabe1}, \cite{watanabe2}. It was proved that the Bayesian generalization 
error has asymptotic expansion given by
\begin{eqnarray*}
\mathbb{E} [ G_n ] = L(w_0) + \frac {\lambda} {n} + o \left ( \frac {1} {n} \right ),
\end{eqnarray*}
where $L(w)$ is the log loss function 
\[
L(w) = - \int q(x) \log p(x|w) dx, 
\]
and $\lambda>0$ is the absolute value of the largest pole of the meromorphic function that is
 analytically continued from 
\[
\zeta(z) = \int K(w)^z \varphi (w) dw.
\]
The constant $\lambda$ is called the learning coefficient in statistics or 
the real log canonical threshold (RLCT) in real algebraic geometry.
The concrete values of RLCTs were studied in normal mixtures \cite{Yam1}, 
reduced rank regressions \cite{aoyagi}, and many statistical models \cite{aoyagi2,aoyagi3,aoyagi4}. 
Mathematical derivation using toric modification was introduced \cite{Yam2}. 
Applications of RLCTs to statistics were also widely studied, for example, model selection by the 
marginal likelihood using RCLT was proposed \cite{Drton}, \cite{2013}, 
and exchange probability  in Markov chain Monte Carlo was analyzed \cite{Nagata}.

In this paper, we study the generalization error of a Poisson mixture and the simplex 
Vandermonde Matrix type singularities, and clarify RLCTs for such statistical models. 
The Poisson distribution is defined by
\begin{eqnarray*}
Po(x| b) = e^{-b} \frac {b^{x}} {{x} !} , \hspace{5mm} x \in \mathbb{Z}_{\ge 0} , \hspace{5mm} b > 0.
\end{eqnarray*}
The M-dimensional Poisson distribution is defined by
\begin{eqnarray*}
{\bf Po(x| b)} =  \prod_{m=1}^M Po(x_{m}| b_{m}), \hspace{5mm} {\bf x} \in \left ( \mathbb{Z}_{\ge 0} \right ) ^ M, \hspace{5mm}
{\bf b} \in \left ( \mathbb{R}_{>0} \right ) ^ M.
\end{eqnarray*}
The M-dimensional Poisson mixture with $H$ components is defined by
\begin{eqnarray*}
p(X=x|{\bf a},{\bf b}) = \sum_{k=1}^H  a_{k}  {\bf Po}({\bf x}| {\bf b}_{k}), \hspace{5mm} \sum_{k=1}^{H} a_{k} = 1, \hspace{5mm}   a_{k} \ge 0.
\end{eqnarray*}
Since the Poisson distribution is defined on $\mathbb{Z}_{\ge 0}$, 
the integration over $x\in \mathbb{R}^M$ is replaced by the summation over $x\in \mathbb{Z}_{\ge 0}$.
Note that $\{a_k\}$ is an element of the $(K-1)$-dimensional symplex, hence singularities
of a mixture models are referred to as simplex Vandermonde type singularities. In this paper we 
give a mathematical derivation for such general cases.

\section{Main Result}

In this section, we explain the main result of this paper. 
We clarify RLCT when $K(w)$ is the mean error function of a Poisson mixture. That is, 
\begin{eqnarray*}
K(w) := - \sum_{x\in \mathbb{Z}_{\ge 0}^M} \left ( \sum_{k=1}^{H^*}  a_{k}^*  {\bf Po}(x| {\bf b}_{k}^*) \right ){\rm log}  \frac { \sum_{k=1}^H  a_{k}  {\bf Po}(x| {\bf b}_{k})} { \sum_{k=1}^{H^*}  a_{k}^*  {\bf Po}(x| {\bf b}_{k}^*)} ,
\end{eqnarray*}
where $a_{k}^*, {\bf b}_{k}^*$ are constant values. 
This function $K(w)$ is equal to the Kullback Leibler divergence from the Poisson mixture with
$H^*=r$ components and that with $H$ components. Also we assume that $\varphi(w)>0$ for 
an arbitrary $w\in W$ which satisties $p(x|w_0)=p(x|w)$ $(\forall x)$. 
Then we prove that
\begin{eqnarray*}
\mathbb{E} [ G_n ] &=&L(w_0) + \frac {3r + H - 2} {4n}  + o \left ( \frac {1} {n} \right )  \ \ \mbox{if dimension}=1,
\\
\mathbb{E} [ G_n ] &=& L(w_0) + \frac {Mr+H-1} {2n}  + o \left ( \frac {1} {n} \right )\ \ \mbox{if dimension}>1.
\end{eqnarray*}
In 2019, Aoyagi derived RLCT of a Poisson mixture model when $H^*=1$\cite{aoyagi4}, however RLCTs of $H^*>1$ have been left unknown. The main result of this paper is to clarify  RLCT of a Poisson mixture
for cases $H^*>1$.
In 2018, we derived an upper bound of the Poisson mixture \cite{sato1}. We prove that it is also equal to the
lower bound, resulting that it is the RLCT. 
\vskip3mm\noindent
{\bf Remark}. Since the free energy satisfies $\mathbb{E}[G_n]=\mathbb{E}[F_{n+1}]-\mathbb{E}[F_n]$,
it follows that
\[
\mathbb{E}[F_n]=nL(w_0)+\lambda\log n +o(\log n).
\]
This result is useful to study singular BIC (sBIC) and WBIC  \cite{Drton}, \cite{2013}.

\section{Preparation \& Previous Research}
\subsection{Polynomial properties}
Let $\mathcal{C}_{r}^{H}$ be the coefficient of $t^{H-r}$ of a polynomial $\prod_{i=1}^H (t+b_i)$.
\begin{ex.}
\begin{eqnarray*}
(t+ b_1) (t+b_2) &=& t^2 + (b_1+b_2) t + b_1 b_2\\
&=& \mathcal{C}_{0}^{2} t^2 + \mathcal{C}_{1}^{2} t + \mathcal{C}_{2}^{2},
\end{eqnarray*}
where $\mathcal{C}_{0}^{2}  := 1,\ \mathcal{C}_{1}^{2}  := b_1+b_2,\ \mathcal{C}_{2}^{2}  := b_1 b_2$.
\end{ex.}
By this notation, 
\begin{eqnarray*}
\prod_{i=1}^H (t+b_i) = \sum_{r=0}^H \mathcal{C}_{r}^{H} t^{H-r}.
\end{eqnarray*}
\begin{th.}
\label{th1}
If $n>H$,
\begin{eqnarray*}
\sum_{r=0}^H (-1)^r \mathcal{C}_{r}^{H} \left ( \sum_{i=1}^H a_i b_i^{n-r} \right )=0.
\end{eqnarray*}
\end{th.}
\begin{proof}
\begin{eqnarray*}
&\sum_{r=0}^H (-1)^r \mathcal{C}_{r}^{H} \left ( \sum_{i=1}^H a_i b_i^{n-r} \right )\\
& = \sum_{i=1}^H a_i b_i^{n-H} \sum_{r=0}^H  \mathcal{C}_{r}^{H} (-1)^r b_i^{H-r}.
\end{eqnarray*}
Therefore, if we can prove $\sum_{r=0}^H  \mathcal{C}_{r}^{H} (-1)^r b_i^{H-r} = 0$, then Theorem \ref{th1} is shoven by the following Lemma \ref{lm1}.\\
\end{proof}
\begin{lm.}
\label{lm1}
$\sum_{r=0}^H  \mathcal{C}_{r}^{H} (-1)^r b_i^{H-r} = 0$.
\end{lm.}
\begin{proof}
By definition of $\mathcal{C}_{r}^{H}$,
\begin{eqnarray*}
\prod_{i=1}^H (t+b_i) = \sum_{r=0}^H \mathcal{C}_{r}^{H} t^{H-r}.
\end{eqnarray*}
Assign $-b_j$ to $t$,
\begin{eqnarray*}
\prod_{i=1}^H (-b_j+b_i) = \sum_{r=0}^H \mathcal{C}_{r}^{H} (-1)^{H-r} b_j^{H-r}.
\end{eqnarray*}
Since $(-b_j+b_j)=0$, $\prod_{i=1}^H (-b_j+b_i) =0$. Thus,
\begin{eqnarray*}
\sum_{r=0}^H \mathcal{C}_{r}^{H} (-1)^{H-r} b_j^{H-r} =0.
\end{eqnarray*}
Finally, by multiplying the both sides of the equation by $(-1)$ to match the sign of Lemma \ref{lm1},  Lemma \ref{lm1} is proved.
\end{proof}
\begin{lm.}
\label{lm2}
$\sum_{i=1}^H a_i b_i^{n} = \sum_{i=1}^H \mathcal{F}_i^{(n)} \left ( \sum_{j=1}^H  a_j b_j^{i} \right ) $,\\
where \[
\mathcal{F}_i^{(n)} := \begin{cases}
\delta_{n, i} & n \le H\\
\sum_{r=1}^H (-1)^{r+1} \mathcal{C}_{r}^{H}  \mathcal{F}_i^{(n-r)} & {\rm Otherwise} 
\end{cases}\ \ \ .
\]
\end{lm.}
\begin{proof}\ \\
The case $n \le H$, it is trivial. \\
The case $n > H$, by Theorem \ref{th1} and assumption,
\begin{eqnarray*}
\sum_{i=1}^H a_i b_i^{n} &=& \sum_{r=1}^H (-1)^{r+1} \mathcal{C}_{r}^{H} \left ( \sum_{i=1}^H a_i b_i^{n-r} \right )\\
&=& \sum_{r=1}^H (-1)^{r+1} \mathcal{C}_{r}^{H} \left (  \sum_{i=1}^H \mathcal{F}_i^{(n-r)} \left ( \sum_{j=1}^H  a_j b_j^{i} \right )  \right )\\
&=&  \sum_{i=1}^H  \left ( \sum_{r=1}^H (-1)^{r+1} \mathcal{C}_{r}^{H}  \mathcal{F}_i^{(n-r)}  \right ) \left ( \sum_{j=1}^H  a_j b_j^{i} \right )  .
\end{eqnarray*}
Thus,
\[
\mathcal{F}_i^{(n)} =  \sum_{r=1}^H (-1)^{r+1} \mathcal{C}_{r}^{H}  \mathcal{F}_i^{(n-r)} .
\]
\end{proof}
Next, we prove multidimensional version of  Lemma \ref{lm2}. Before proving that, introduce Multi-index notation to simplify formulas as below.
\begin{eqnarray*}
{\bf b}^{\bf r} = \prod_{m=1}^M b_m^{r_m},  {\bf b} \in \mathbb{R}^M, {\bf r} \in \mathbb{Z}_{\ge 0}^{M}.
\end{eqnarray*}
Let $\mathcal{C}_{r}^{H} (b_{\cdot m})$ be the coefficient of $t^{H-r}$ of $\prod_{i=1}^H (t+b_{im})$.
Define $\mathcal{F}_i^{(n)} (b_{\cdot m})$ by
\[
\mathcal{F}_i^{(n)} (b_{\cdot m}) := \begin{cases}
\delta_{n, i} & n \le H\\
\sum_{r=1}^H (-1)^{r+1} \mathcal{C}_{r}^{H} (b_{\cdot m}) \mathcal{F}_i^{(n-r)} (b_{\cdot m})& {\rm Otherwise} 
\end{cases},
\]
and define $\mathcal{F}_{\bf r}^{(\bf n)}$ using $\mathcal{F}_i^{(n)} (b_{\cdot m})$ by
\[
\mathcal{F}_{\bf r}^{(\bf n)} := \prod_{m=1}^M \mathcal{F}_{ r_m}^{(n_m)} (b_{\cdot m}).
\]

\begin{lm.}
\label{lm3}
$\sum_{i=1}^H a_i {\bf b}_i^{\bf n} = \sum_{{\bf r} \in [1:H]^M} \mathcal{F}_{\bf r}^{(\bf n)} \left (  \sum_{i=1}^H a_i {\bf b}_i^{\bf r} \right ) $.\\
In other words, $\sum_{i=1}^H a_i \prod_{m=1}^M b_{im}^{n_m} = \sum_{{\bf r} \in [1:H]^M} \mathcal{F}_{\bf r}^{(\bf n)} \left (  \sum_{i=1}^H a_i \prod_{m=1}^M b_{im}^{r_m} \right ) $.\\
\end{lm.}
\begin{proof}\ \\
We prove this lemma by mathematical induction.
In the case that the dimension of $b$ equals to $1$, it follows from Lemma \ref{lm2}.
Assume that if the dimension of $b$ equals to $M-1$, Lemma \ref{lm3} holds. Then, if the dimension of $b$ equals to $M$, we prove Lemma \ref{lm3} as follows.
Define $ c_{k}^{({\bf n})} := a_{k} \prod_{m=1}^{M-1} b_{km}^{n_m}$ for all $k\in [1:H]$. By using Lemma \ref{lm2},
\begin{eqnarray*}
\sum_{i=1}^H a_i {\bf b}_i^{\bf n} &=& \sum_{i=1}^H c_i^{({\bf n})} b_{i M}^{n_M} \\
& =& \sum_{i=1}^H \mathcal{F}_i^{(n_M)} (b_{\cdot M}) \left ( \sum_{j=1}^H  c_j^{({\bf n})} b_{jM}^{i} \right ) \\
& =& \sum_{i=1}^H \mathcal{F}_i^{(n_M)} (b_{\cdot M}) \left ( \sum_{j=1}^H  a_{j} b_{jM}^{i} \prod_{m=1}^{M-1} b_{jm}^{n_m}  \right ) .\\
\end{eqnarray*}
Define $ d_{j}^{(i)} := a_{j} b_{jM}^i$ for all $i, j\in [1:H]$.
\begin{eqnarray*}
\sum_{i=1}^H a_i {\bf b}_i^{\bf n} &=& \sum_{i=1}^H \mathcal{F}_i^{(n_M)} (b_{\cdot M}) \left ( \sum_{j=1}^H  d_j^{(i)} \prod_{m=1}^{M-1} b_{jm}^{n_m}  \right ) \\
& =& \sum_{i=1}^H \mathcal{F}_i^{(n_M)} (b_{\cdot M}) \left (\sum_{{\bf r} \in [1:H]^{M-1}} \left ( \prod_{m=1}^{M-1}  \mathcal{F}_{r_m}^{(n_m)} (b_{\cdot m}) \right ) \left (  \sum_{k=1}^H d_k^{(i)} \prod_{m=1}^{M-1} b_{km}^{r_m} \right ) \right ) \\
& =& \sum_{i=1}^H\sum_{{\bf r} \in [1:H]^{M-1}} \left ( \prod_{m=1}^{M}  \mathcal{F}_{r_m}^{(n_m)} (b_{\cdot m}) \right ) \left (  \sum_{k=1}^H  a_{k} b_{kM}^i \prod_{m=1}^{M-1} b_{km}^{r_m} \right ) \\
& =& \sum_{{\bf r} \in [1:H]^{M}} \mathcal{F}_{\bf r}^{(\bf n)}  \left (  \sum_{k=1}^H  a_{k}  {\bf b}_{k}^{r} \right ) , \\
\end{eqnarray*}
where we used the assumption in the second equation.
\end{proof}
\subsection{Learning coefficient properties}
[This section is based on the book "Algebraic Geometry and Statistical Learning Theory"\cite{hon2}]\\
Let $\lambda(K,\varphi)$ denote the learning coefficient(RLCT) of $K(w), \varphi (w)$ where $K(w)$ is an analytic function. It is shown that LRCT satisfies the following equations and inequalities. (The last property is original, the other is written in the book or his paper.)
\begin{description}
\item[Sum]
If $w=(w_1, w_2), K(w)=K_1 (w_1) + K_2 (w_2), \varphi (w) = \varphi_1 (w_1) \varphi_2 (w_2)$,
\begin{eqnarray*}
\lambda (K, \varphi) = \lambda ( K_1, \varphi_1) + \lambda (K_2, \varphi_2).
\end{eqnarray*}
\item[Product]
If $w=(w_1, w_2), K(w)=K_1 (w_1)  K_2 (w_2), \varphi (w) = \varphi_1 (w_1) \varphi_2 (w_2)$,
\begin{eqnarray*}
\lambda (K, \varphi) = \min \left \{ \lambda ( K_1, \varphi_1) , \lambda (K_2, \varphi_2) \right \} .
\end{eqnarray*}
\item[Inequality]
If $K_1 (w) \le K_2 (w)$ and $ \varphi_1 (w) \ge \varphi_2 (w)$,
\begin{eqnarray*}
\lambda (K_1, \varphi_1) \le \lambda (K_2, \varphi_2).
\end{eqnarray*}
\item[Inequality: constant factor]
If there exists $c_1, c_2 \in \mathbb{R}_{>0}$ such that $c_1 K_1(w) \le K_2(w) \le c_2 K_1 (w)$,
\begin{eqnarray*}
\lambda (K_1, \varphi) = \lambda (K_2, \varphi).
\end{eqnarray*}
\item[Bounded function]
Let $\{f_l (w)\}_{l\in[1: L]}$ and $\{g_p (w) \}_{p\in [1:R]}$ be sets of analytic functions and define $K_1 (w) := \sum_{l=1}^{L} f_l (w)^2$ and $K_2(w): = \sum_{p=1}^{R} g_p (w)^2$.
If there exist bounded functions $\{h_{lp} (w)\}_{l \in [1:L], p \in [1:R]}$ on $W$ such that $f_l(w) = \sum_{p=1}^{R}  g_p (w) h_{lp}(w)$,
\[\lambda (K_1, \varphi) \le \lambda (K_2, \varphi ). \]
\begin{proof}\ \\
By using Cauchy-Schwarz inequality,
\begin{eqnarray*}
K_1 (w) &=&  \sum_{l=1}^{L} f_l (w)^2\\
&=&  \sum_{l=1}^{L} \left ( \sum_{p=1}^{R}  g_p (w) h_{lp}(w) \right )^2\\
&\le& \sum_{l=1}^{L} \left ( \sum_{p=1}^{R}    g_p (w) ^2 \right ) \left ( \sum_{p=1}^{R}  h_{lp}(w)^2 \right )\\
&=& \left (  \sum_{l=1}^{L} \sum_{p=1}^{R}  h_{lp}(w)^2 \right ) \left ( \sum_{p=1}^{R}    g_p (w) ^2 \right ) \\
&\le& M K_2(w),
\end{eqnarray*}
where $M$ is contant. Then, by applying properties of inequality, property of Bounded function holds.
\end{proof}
\item[Ideal Invariance]
If the ideal generated from $\{f_l (w)\}_{l=[1: L]}$ and the ideal generated from $\{g_m(w)\}_{m\in [1: M]}$ are equivalent in a ring of convergent power series $\mathbb{R}\langle W \rangle $ and define
\begin{eqnarray*}
K_1 (w) := \sum_{l=1}^{L} f_l (w)^2 ,K_2 (w) := \sum_{m=1}^M g_m(w)^2,
\end{eqnarray*}
then
\begin{eqnarray*}
\lambda (K_1, \varphi) = \lambda (K_2, \varphi).
\end{eqnarray*}
\end{description}
\subsubsection{Notation}
We introduce new notation$=_{RLCT}$ to formula and ideal.\\
Let $f,g$ be analytic functions.\\
An equivalence relation $=_{RLCT}$ is defined by
\begin{eqnarray*}
& \ &f(w)  =_{RLCT} g(w) \\
&\iff & \exists C_1, C_2 > 0, \ \  \forall w \in W \ \left ( C_1 f(w)^2 \le g(w)^2 \le C_2 f(w)^2 \right ) .
\end{eqnarray*}
For ideals, $=_{RLCT}$ is defined by
\begin{eqnarray*}
& \ &\langle \{f_l (w)\}_{l=[1: L]}\rangle =_{RLCT} \langle \{g_m(w)\}_{m\in [1: M]} \rangle\\
&\iff & \sum_{l=1}^{L} f_l (w)^2=_{RLCT} \sum_{m=1}^M g_m(w)^2.
\end{eqnarray*}
Moreover,  we introduce a new notation $\le_{RLCT}$,
\begin{eqnarray*}
& \ &f(w)  \le_{RLCT} g(w) \\
&\iff & \exists C_1 > 0, \ \  \forall w \in W \ \left ( f(w)^2 \le C_1 g(w)^2 \right ) .
\end{eqnarray*}
\begin{eqnarray*}
& \ &\langle \{f_l (w)\}_{l=[1: L]}\rangle \le_{RLCT} \langle \{g_m(w)\}_{m\in [1: M]} \rangle\\
&\iff &  \sum_{l=1}^{L} f_l (w)^2  \le_{RLCT} \sum_{m=1}^M g_m(w)^2.
\end{eqnarray*}
The definition above is useful to express large/small relationship on the whole parameter space $W$ but sometimes we want to express large/small relationship on a local parameter region, for example, when calculating RLCT,  we will restrict parameter to a neighborhood of a certain point $w^* \in W$. In such cases, we will use the word "Locally" or "on a neighborhood of $w^*$" to express the restriction of parameter region. \\

The last useful property of RLCT is given by the following Lemma.\\
\begin{lm.}
Let $\langle \{f_l (w)\}_{l=[1: L]} \rangle$ be the ideal generated from $\{f_l (w)\}_{l=[1: L]}$ in a ring of convergent power series $\mathbb{R}\langle W \rangle $.
If there exists $f_{L+1}(w) \in\mathbb{R}\langle W \rangle $ such that
\begin{eqnarray*}
f_{L+1}(w)^2 \le \sum_{l=1}^L f_l (w)^2 \ \ \ \ \ on\ W,
\end{eqnarray*}
then 
\[
\langle \{f_l (w)\}_{l=[1: L]} \rangle =_{RLCT} \langle \{f_l (w)\}_{l=[1: L+1]} \rangle.
\]
\end{lm.}
\begin{proof}
\ \\
$\langle \{f_l (w)\}_{l=[1: L]}  \rangle \le_{RLCT} \langle \{f_l (w)\}_{l=[1: L+1]} \rangle$ follows from $\sum_{l=1}^L f_l(w)^2 \le \sum_{l=1}^{L+1} f_l(w)^2$ and property of inequality.And
\begin{eqnarray*}
\sum_{l=1}^{L+1} f_l(w)^2 \le \sum_{l=1}^{L} f_l(w)^2 + \sum_{l=1}^L f_l (w)^2 =2  \sum_{l=1}^{L} f_l(w)^2 \ \ \ \ \ on\ W.
\end{eqnarray*}
Therefore, apply property of inequality, $ \langle \{f_l (w)\}_{l=[1: L]}  \rangle \le_{RLCT} \langle \{f_l (w)\}_{l=[1: L+1]}  \rangle \le_{RLCT} \langle \{f_l (w)\}_{l=[1: L]}  \rangle$.
\end{proof}
\subsection{Poisson Mixture Properties}
The mean error function of a Poisson mixture is not analytic function of the parameter. However, in 2004, asymptotic expansion of the stochastic complexity of non-analytic learning machines was established by Watanabe\cite{noana}. 
\begin{cd.} \label{cd:1}
The parameter space $W$ is $W \subset \mathbb{R}^d$ and compact. Using the $C^{\infty}$ function $\varphi_0 (w)$ that is always positive on $W$ and the analytic function $\varphi_1 (w)$ that is nonnegative, the prior distribution $\varphi (w)$ is defined as $\varphi (w) = \varphi_0 (w) \varphi_1 (w) $.
\end{cd.}
\begin{df.}
Define $S(t) := \begin{cases} \frac {- {\rm log} t + t -1} {(t-1)^2} & (t>0 \& t \not = 1)\\ 1 &(t = 1)\end{cases} $ . This function is a positive and analytic function.
\end{df.}
\begin{cd.} \label{cd:2}
There exists measurable function $M(x)$ on $\mathbb{R}^m$ such that
\[
\sup_{w \in W} S \left ( \frac {p(x|w)} {q(x)}\right ) \le M(x).
\]
And define Lebesgue-Stieltjes measure $\mu$ by $\mu (a,b] := \int_{(a, b]} M(x) q(x) dx$. We can define an inner product as below.
\[
(u, v) := \int u v d \mu = \int u(x) v(x) M(x) q(x) dx.
\]
 Assume $\Phi$ is analytic of $w$ in Hilbert space $L^2 (X, \mu)$.
\[
\Phi : W \ni w \mapsto \frac {p(x|w)} {q(x)} - 1 \in L^2 (X, \mu).
\]
\end{cd.}
\begin{th.}
\label{th:prev:02}
Suppose Condition \ref{cd:1} and Condition \ref{cd:2} and define
\[
K(w) := \int  \left [ \frac {p(x|w)} {q(x)} - 1 \right ]^2 M(x) q(x) dx.
\]
The zeta function of $K(w)$ is defined by
\[
\zeta (z) := \int {K(w)}^z \phi(w) dw .
\]
Let $-\lambda$ be the largest pole of $\zeta (z) $ and $m$ be its order.  The asymptotic expansion of the free energy $F_n$ is 
\[
F_n = \lambda {\log} n - (m-1) {\rm log} {\rm log} n + R + o_p(1),
\]
where $R$ is a random variable.
$\int_C  \left [ p(x|w)  - q(x) \right ]^2 dx$ has the same pole of $K(w)$ where $\exists C$ in a domain of $x$.
\end{th.}
\subsubsection{Check that Poisson mixture model satisfies Condition \ref{cd:1} and Condition \ref{cd:2}.}
\begin{lm.}
\label{112301}
Let $p(x|w)$,$q(x)$ be Poisson mixtures. There exist constant values $A_0, A_1, B_0, B_1>1$ such that
\begin{eqnarray*}
	\frac {1} {A_0} \prod_{m=1}^{M} \left ( \frac {1} {A_1} \right ) ^{x_m} \le \frac {p(x|w)} {q(x)} \le  B_0 \prod_{m=1}^{M} B_1^{x_m}
\end{eqnarray*}
for all $w \in W$, for all $x \in \mathbb{R}^N$.
\end{lm.}
\begin{proof}
\begin{eqnarray*}
\frac {p(x|w)} {q(x)} &=& \frac {\sum_{k=1}^{H} a_k \prod_{m=1}^{M} e^{-b_{km}} \frac {b_{km}^{x_m}} {x_m!}} {\sum_{k=1}^{H^*} {a^*}_k \prod_{m=1}^{M} e^{-{b^*}_{km}} \frac {{b^*}_{km}^{x_m}} {x_m!}} 
\le  \frac {\sum_{k=1}^{H} a_k \prod_{m=1}^{M} e^{-b_{km}} b_{km}^{x_m}} {{a^*}_1 \prod_{m=1}^{M} e^{-{b^*}_{1m}} {b^*}_{1m}^{x_m} } \\
&=& \frac {1} {{a^*}_1} \sum_{k=1}^{H} a_k \prod_{m=1}^{M} e^{{b^*}_{1m}-b_{km}} \left ( \frac {b_{km}} {{b^*}_{1m}} \right )^{x_m}  \\
& \le & \frac { \max_{m \in {1,2,...,M}} \{ e^{M {b^*}_{1m}} \} } {{a^*}_1} \sum_{k=1}^{H} a_k \prod_{m=1}^{M}  \left ( \frac {b_{km}} {{b^*}_{1m}} \right )^{x_m} . \\
\end{eqnarray*}
Let $B_0=\frac {\max_{m \in {1,2,...,M}} \{ e^{M {b^*}_{1m}} \}} {{a^*}_1}$. Since the parameter space $W$ is a bounded closed set, there exsits  $B_1 \ge 1$ such that $\frac {b_{km}} {{b^*}_{1m}} \le B_1$ for all $k,m$.
As a result,
\begin{eqnarray*}
\frac {p(x|w)} {q(x)} & \le & B_0 \sum_{k=1}^{H} a_k \prod_{m=1}^{M}  {B_1}^{x_m}  B_0 \prod_{m=1}^{M}  {B_1}^{x_m} .
\end{eqnarray*}
We prove the right side inequality of lemma and we can prove the left side in a similar way forcusing on inversion. 
\end{proof}
\begin{lm.}
\label{112302}
Let $p(x|w)$,$q(x)$ be Poisson mixtures.
There exist constant values $C_0, C_1>0$ such that 
\begin{eqnarray*}
S \left ( \frac {p(x|w)} {q(x)} \right ) \le C_0 \sum_{k=1}^{H} x_k +C_1
\end{eqnarray*}
for all $w\in W$, for all $x\in \mathbb{R}^N$.
\end{lm.}
\begin{proof}
Because of the definition of $S(t)$, $S(t)>0(0<t<\infty)$ is a continuous function and $S(t)$ diverge to the positive infinity when $t \to 0$. Thus, there exist $D_0,D_1>0$ such that
\begin{eqnarray*}
S(t) \le \max \left ( -D_0 {\rm log} t,\ D_1 \right ).
\end{eqnarray*}
From the left side of Lemma \ref{112301},
\begin{eqnarray*}
\max \left ( -D_0 {\rm log} \frac {p(x|w)} {q(x)} ,\ D_1 \right )&\le& \max [ D_0 {\rm log} A_0 + D_0 \left ( \sum_{m=1}^M x_m \right ) {\rm log} A_1, D_1 ] \\
& \le & \left ( D_0 {\rm log} A_1 \right ) \left ( \sum_{m=1}^M x_m \right ) +D_0 {\rm log} A_0 + D_1 .
\end{eqnarray*}
Therefore, Lemma \ref{112302} follows from the above inequality and $ D_0 {\rm log} A_1 ,D_0 {\rm log} A_0 + D_1>0$.
\end{proof}
\begin{th.}
Let $p(x|w)$,$q(x)$ be Poisson mixtures.
Using $C_0, C_1>0$ of Lemma \ref{112302}, define $M(x)$ by
\begin{eqnarray*}
M(x) = C_0 \sum_{m=1}^M x_m + C_1.
\end{eqnarray*}
Then $\frac {p(x|w)} {q(x)}-1$ satisfies Condition \ref{cd:2}.
\end{th.}
\begin{proof}
Let  $u(x)$ be as below.
\begin{eqnarray*}
u(x) &:=& \frac {p(x|w)} {q(x)} -1 =  \frac {\sum_{k=1}^{H} a_k \prod_{m=1}^{M} e^{-b_{km}} b_{km}^{x_m}} {\sum_{k=1}^{H^*} {a^*}_k \prod_{m=1}^{M} e^{-{b^*}_{km}} {b^*}_{km}^{x_m} } - 1 .
\end{eqnarray*}
Taylor expansion of $b_{km}^{x_m}$ is
\begin{eqnarray*}
{b_{km}}^{x_m} = \sum_{j=0}^{\infty} \frac {1} {j!} \left ( x_m {\rm log} b_{km} \right )^{j}.
\end{eqnarray*}
Taylor expansion of $u(x)$ to $J$-th term is
\begin{eqnarray*}
u_J (x,w) := \sum_{j=0}^{J} \frac {1} {j!} \frac {\sum_{k=1}^{H} a_k \prod_{m=1}^{M} e^{-b_{km}} \left ( x_m {\rm log} b_{km}  \right )^{j}} {\sum_{k=1}^{H^*} {a^*}_k \prod_{m=1}^{M} e^{-{b^*}_{km}} {b^*}_{km}^{x_m} } - 1.
\end{eqnarray*}
The function ${\rm log} y $ is an analytic function on $y>0$ so $u_J (x,w)$ is an analytic function of the parameter $w(=({\bf a}, {\bf b}))$. If we fix the parameter $w$, then $u_J(\cdot,w)$ is an element of Hilbert space and then  $u_J(x,w)$ is an analytic function whose co-domain is Hilbert space. If we prove that  $u_J (x,w)$ uniformly converges to $u(x,w)$ in Hilbert space($J \to \infty$), Condition \ref{cd:2} follows from that, if an analytic function $f_n$ uniformly converges to $f$, then $f$ is also an analytic function.  We prove $u_J (x,w) \to u(x,w)$($J \to \infty$) uniformly in Hilbert space.
\begin{eqnarray*}
T_J & \equiv & \sup_{w \in W} || u(w)- u_J (w) ||^2  \\
&=& \sup_{w \in W} \sum_{x=0}^{\infty} q(x) M(x) \left [   \sum_{j>J} \frac {1} {j!} \frac {\sum_{k=1}^{H} a_k \prod_{m=1}^{M} e^{-b_{km}} \left ( x_m {\rm log} b_{km}  \right )^{j}} {\sum_{k=1}^{H^*} {a^*}_k \prod_{m=1}^{M} e^{-{b^*}_{km}} {b^*}_{km}^{x_m} } \right ]^2 \\
& \le & \sup_{w \in W} \sum_{x=0}^{\infty} q(x) M(x) \left [   \sum_{j>J} \frac {1} {j! {a^*}_1 } \sum_{k=1}^{H} a_k \prod_{m=1}^{M} e^{{b^*}_{1m}-b_{km}} \frac { \left ( x_m {\rm log} b_{km}  \right )^{j} } { {b^*}_{1m}^{x_m}} \right ]^2    \\
& \le &  \sum_{x=0}^{\infty} q(x) M(x) \left [   \sum_{j>J} \frac {1} {j!} B_0 \prod_{m=1}^{M} B_1^{x_m}  {x_m}^j \right ]^2,
\end{eqnarray*}
where $B_0,B_1$  are same values the right side of Lemma \ref{112301}.
Therefore
\begin{eqnarray*}
\lim_{J \to \infty} T_{J} & \le & \lim_{J \to \infty} \sum_{x=0}^{\infty} q(x) M(x) \left [   \sum_{j>J} \frac {1} {j!} B_0 \prod_{m=1}^{M} B_1^{x_m}  {x_m}^j \right ]^2. \\
\end{eqnarray*}
If the dominated convergence theorem can be employed,
\begin{eqnarray*}
\lim_{J \to \infty} T_{J} & \le &\sum_{x=0}^{\infty}  \lim_{J \to \infty}  q(x) M(x) \left [   \sum_{j>J} \frac {1} {j!} B_0 \prod_{m=1}^{M} B_1^{x_m}  {x_m}^j \right ]^2. \\
\end{eqnarray*}
The type of the rightmost factor is similar to Poisson mixture whose parameter is $\prod_{m=1}^M x_m$ so $\sum_{j=0}^{\infty} Po \left (j | \prod_{m=1}^M x_m \right ) = 1$
\begin{eqnarray*}
\lim_{J \to \infty} T_{J} \le \sum_{x=0}^{\infty}  0 
= 0.
\end{eqnarray*}
The dominated convergence theorem can be employed because there exists a Lebesgue integrable function $g(x)$ which dominates $  \left [   \sum_{j>J} \frac {1} {j!} B_0 \prod_{m=1}^{M} B_1^{x_m}  {x_m}^j \right ]^2$.
\begin{eqnarray*}
\left [   \sum_{j>J} \frac {1} {j!} B_0 \prod_{m=1}^{M} B_1^{x_m}  {x_m}^j \right ]^2  &=& \left [  B_0 \prod_{m=1}^{M} e^{x_m} B_1^{x_m}    \sum_{j>J} e^{-x_m} \frac {{x_m}^j} {j!}  \right ]^2 \\
&\le&  \left [ B_0 \prod_{m=1}^{M} e^{x_m} B_1^{x_m}    \sum_{j=0}^{\infty} e^{-x_m} \frac {{x_m}^j} {j!}  \right ]^2 \\
&=& \left [ B_0 \prod_{m=1}^{M} \left ( e B_1 \right )^{x_m}  \right ]^2 = B_0^2 \prod_{m=1}^{M} \left ( e^2 B_1^2 \right )^{x_m}  . 
\end{eqnarray*}
Let $U_0 := B_0^2,\ U_1 :=e^2 B_1^2 $, $B_0,\ B_1,\ e$ abe constants so $U_0,\ U_1$ are also constants.
\begin{eqnarray*}
B_0^2 \prod_{m=1}^{M} \left ( e^2 B_1^2 \right )^{x_m} &=& U_0 \prod_{m=1}^{M} U_1 ^{x_m} .
\end{eqnarray*}
Define $g(x):= U_0 \prod_{m=1}^{M} U_1 ^{x_m} $. From the definition, it is clear that $ \left [   \sum_{j>J} \frac {1} {j!} B_0 \prod_{m=1}^{M} B_1^{x_m}  {x_m}^j \right ]^2\le g(x)$. Next, we prove $g(x)$ has a finite sum.
\begin{eqnarray*}
\sum_{x=0}^{\infty} g(x) &=& \sum_{x=0}^{\infty} q(x)M(x) U_0 \prod_{m=1}^{M}  U_1^{x_m}\\
 &=& \sum_{x=0}^{\infty} \left ( \sum_{k=1}^{H^*} {a^*}_k \prod_{m=1}^{M} e^{-{b^*}_{km}} \frac {{b^*}_{km}^{x_m}} {x_m!} \right ) M(x)   U_0 \prod_{m=1}^{M}  U_1 ^{x_m}  \\
 &=& \sum_{x=0}^{\infty} \left ( \sum_{k=1}^{H^*} {a^*}_k \prod_{m=1}^{M} e^{-{b^*}_{km}} \frac { \left ({b^*}_{km} U_1 \right )^{x_m} } {x_m !}\right ) M(x)   U_0   \\
 &=& \sum_{x=0}^{\infty} \left ( \sum_{k=1}^{H^*} {a^*}_k \prod_{m=1}^{M}  e^{{b^*}_{km} (U_1-1) }  Po( x_m |{b^*}_{km} U_1 ) \right ) M(x)   U_0   \\
 &=&  U_0 \sum_{k=1}^{H^*} \sum_{x=0}^{\infty}M(x)   {a^*}_k  e^{(U_1-1) \sum_{m=1}^M {b^*}_{km} } \prod_{m=1}^{M}  Po( x_m |{b^*}_{km} U_1 )     \\
 &=& U_0 \sum_{k=1}^{H^*}  {a^*}_k  e^{(U_1-1) \sum_{m=1}^M {b^*}_{km} }\sum_{x=0}^{\infty}M(x)   \prod_{m=1}^{M} Po( x_m |{b^*}_{km} U_1 )     \\
 &=& U_0 \sum_{k=1}^{H^*}  {a^*}_k  e^{(U_1-1) \sum_{m=1}^M {b^*}_{km} } \sum_{x=0}^{\infty} \left ( C_0 \sum_{m=1}^M x_m + C_1 \right ) \prod_{m=1}^{M}  Po( x_m |{b^*}_{km} U_1 )     \\
 &=& U_0 \sum_{k=1}^{H^*}  {a^*}_k  e^{(U_1-1) \sum_{m=1}^M {b^*}_{km} } \left ( C_0 U_1 \sum_{m=1}^M {b^*}_{km}  + C_1 \right ) .
\end{eqnarray*}
Since $C_0,\ C_1,\ U_0,\ U_1,\ a^*,\ b^*$ are constants, $g(x)$ has a finite sum.\\
\end{proof}
For deriving learning coefficient (RLCT) of a Poisson mixture, learning coefficient (RLCT) of 
\begin{eqnarray*}
\int_C  \left [ p(x|w)  - q(x) \right ]^2 dx
\end{eqnarray*}
can be used instead of  learning coefficient (RLCT) of 
\begin{eqnarray*}
- \sum_{x\in \mathbb{Z}_{\ge 0}^M} \left ( \sum_{k=1}^{H^*}  a_{k}^*  {\bf Po}(x| {\bf b}_{k}^*) \right ){\rm log}  \frac { \sum_{k=1}^H  a_{k}  {\bf Po}(x| {\bf b}_{k})} { \sum_{k=1}^{H^*}  a_{k}^*  {\bf Po}(x| {\bf b}_{k}^*)} .
\end{eqnarray*}
We write $K(w) := \int_C  \left [ p(x|w)  - q(x) \right ]^2 dx$ hereafter.
\section{Ideal of A Poisson Mixture}
In this section, we study the properties of ideal of a Poisson mixture.
The ideal $\left \langle \left ( p(x|w)-q(x) \right )^2 \right \rangle_{x \in \mathbb{Z}_{\ge 0}^M}$ on $\mathbb{R}$-value convergent power series is given by
\begin{eqnarray*}
 \left ( p(x|w)-q(x) \right )^2 &=& \left ( \sum_{k=1}^H  a_{k}  {\bf Po}(x| b_{k}) -\sum_{k=1}^{H^*}  a_{k}^*  {\bf Po}(x| b_{k}^*) \right )^2\\
 &=& \left ( \sum_{k=1}^H  a_{k}  \prod_{m=1}^M Po(x_{m}| b_{km}) -\sum_{k=1}^{H^*}  a_{k}^*  \prod_{m=1}^M Po(x_{m}| b_{km}^*) \right )^2\\
 &=& \left ( \sum_{k=1}^H  a_{k}  \prod_{m=1}^M  e^{-b_{km}} \frac {b_{km}^{x_m}} {{x_m} !} -\sum_{k=1}^{H^*}  a_{k}^*  \prod_{m=1}^M e^{-b_{km}^*} \frac { {b_{km}^*}^{x_m}} {{x_m} !} \right )^2.\\
\end{eqnarray*}
We show this ideal can be written simply using other generating set.
\subsection{1 dimensional case}
\begin{eqnarray*}
K(w) = \sum_{x=0}^\infty \left ( p(x|w)-q(x) \right )^2 &=&\sum_{x=0}^\infty \left ( \sum_{k=1}^H  a_{k}    e^{-b_k} \frac {b_k^{x}} {x !} -\sum_{k=1}^{H^*}  a_{k}^*  e^{-b_k^*} \frac { {b_k^*}^{x}} {x !} \right )^2.
\end{eqnarray*}
We will transform $K(w)$ to a simple formula in two-phase.
\begin{lm.}
\label{mixed:poisson:lm01}
$K(w)$ is equivalent to
\begin{eqnarray*}
H'(w) = \sum_{x \in [0:H+H^*-1]} \left ( \sum_{k=1}^H a_k e^{-b_k} b_k^x -\sum_{k=1}^{H^*} a_k^* e^{-b_k^*} {b_k^*}^x \right )^2
\end{eqnarray*}
with respect to RLCT.
\end{lm.}
\begin{proof}
\ \\
Firstly, $\lambda(K,\varphi) \le \lambda(H',\varphi)$ is proved by
\[
\sum_{x=0}^{{\bf H+H^*-1}} \left ( \sum_{k=1}^H  a_{k}    e^{-b_k} \frac {b_k^{x}} {x !} -\sum_{k=1}^{H^*}  a_{k}^*  e^{-b_k^*} \frac { {b_k^*}^{x}} {x !} \right )^2
\le \sum_{x=0}^{{\bf  \infty}} \left ( \sum_{k=1}^H  a_{k}    e^{-b_k} \frac {b_k^{x}} {x !} -\sum_{k=1}^{H^*}  a_{k}^*  e^{-b_k^*} \frac { {b_k^*}^{x}} {x !} \right )^2
\]
and RLCT's inequality property. Next, we will prove $\lambda(K,\varphi) \ge \lambda(H',\varphi)$.
Only in this proof, we use new notation as below.
\begin{itemize}
\item $c_k := a_{k} e^{-b_k}$ for all $k\in [1:H]$
\item $ c_{H+k} := a_{k}^* e^{-b_k^*}$ for all $k\in [1:H^*]$
\item $ b_{H+k} := b_k^*$for all $k\in [1:H^*]$
\item $f_n := \sum_{k=1}^{H+H^*} c_k b_k^n$ for all $n \in \mathbb{N}$
\end{itemize}
By using this notation,  $H'(w)$ and $K(w)$ are written simply like $H'(w) = \sum_{n=0}^{H+H^*-1} f_n^2$, $K(w) = \sum_{n=0}^{\infty} \left ( \frac {f_n} {n!}  \right )^2$.\\
Because of Lemma \ref{lm1}, $f_n$ can be written as sum of $H+H^*$ polynomials $f_{n-1},f_{n-2}, \cdots, f_{n-(H+H^*)}$, that is, $f_n = \sum_{i=0}^{H+H^*-1} \mathcal{F}_i^{(n)}  f_{i} $.
Moreover there exists $R$ such that $ |\mathcal{F}_i^{(n)} |<R^n$ for all $i,n$.
The existence of $R$ is proved as below. The definition is
\[\mathcal{F}_i^{(n)} := \begin{cases}
\delta_{n, i} &  n \le {H+H^*-1}\\
\sum_{r=1}^{H+H^*} (-1)^{r+1} \mathcal{C}_{r}^{H+H^*}  \mathcal{F}_i^{(n-r)} & {\rm Otherwise} 
\end{cases}.
\]
Define $R := (H+H') \max_{w \in W}\prod_{i=1}^{H+H^*} (1+|b_i|)  $, then for all $r\in [0:H+H^*]$,
\[
\left | \mathcal{C}_{r}^{H+H^*} \right | \le \sum_{i=0}^{H+H^*} \left | \mathcal{C}_{r}^{H+H^*} \right | \le  \prod_{i=1}^{H+H^*} (1+|b_i|) \le \frac {R} {H+H^*}.
\]
Next, using mathematical induction, we prove $|\mathcal{F}_i^{(n)} |<R^n$.
In the case $n=0,1,2,\cdots, H+H^*-1$, it is trivial because of the definition of $R$.
Assume that $|\mathcal{F}_i^{(k)} |<R^k$ for all $k\in [1:n]$,
\[
\left | \mathcal{F}_i^{(n+1)}  \right | \le \sum_{r=1}^{H+H^*} \left |  \mathcal{C}_{r}^{H+H^*} \right | \cdot \left | \mathcal{F}_i^{(n+1-r)}
\right |\le  \sum_{r=1}^{H+H^*}  \frac {R} {H+H^*}  \cdot  \left | \mathcal{F}_i^{(n+1-r)}
\right |  \overset {\mbox{Assumption}} {\le} \frac {1} {H+H^*} \sum_{r=1}^{H+H^*}  R^{n+2-r}.\]
And $R\ge 1$, $|\mathcal{F}_i^{(n+1)} |<R^{n+1}$.
We proved $ |\mathcal{F}_i^{(n)} |<R^n$ and then we will prove $C H'(w) $ dominates $K(w) = \sum_{n=0}^{\infty} \left ( \frac {f_n} {n!}  \right )^2$ where $C$ is a constant value.
\begin{eqnarray*}
\sum_{n=0}^{\infty} \left ( \frac {f_n} {n!}  \right )^2 &=& \sum_{n=0}^{\infty} \left ( \sum_{i=0}^{H+H^*-1} \frac {  \mathcal{F}_i^{(n)}} {n!}  f_{i}  \right )^2 \\ 
&\le& \sum_{n=0}^{\infty}  \left ( \sum_{i=0}^{H+H^*-1} \left ( \frac {\mathcal{F}_i^{(n)}} {n!} \right )^2 \right )  \left ( \sum_{i=0}^{H+H^*-1} f_{i}^2  \right ) \\ 
&\le&\left ( \sum_{i=0}^{H+H^*-1} f_{i}^2  \right )  \left ( \sum_{i=0}^{H+H^*-1}\sum_{n=0}^{\infty}   \left ( \frac {R^n} {n!} \right )^2 \right ).  \\ 
\end{eqnarray*}
In general, if $a_i \ge 0$, then $\sum_{i=0}^n a_i^2 \le \left ( \sum_{i=0}^n a_i \right ) ^2$.
\begin{eqnarray*}
\sum_{n=0}^{\infty} \left ( \frac {f_n} {n!}  \right )^2 &\le&H'(w) \left ( \sum_{i=0}^{H+H^*-1}   \left (\sum_{n=0}^{\infty} \frac {R^n} {n!} \right )^2 \right )  \\ 
&\le& \left ( H+H^* \right ) e^{2R}H'(w). 
\end{eqnarray*}
Thus, $C H'(w) $ dominates $K(w) = \sum_{n=0}^{\infty} \left ( \frac {f_n} {n!}  \right )^2$.
As a result, $\lambda(K,\varphi) = \lambda(H',\varphi)$.
\end{proof}

\begin{lm.}
$K(w)$ is equivalent to
\begin{eqnarray*}
H(w) = \sum_{x \in [0:H+H^*-1]} \left ( \sum_{k=1}^H a_k  b_k^x -\sum_{k=1}^{H^*} a_k^* {b_k^*}^x \right )^2
\end{eqnarray*}
with respect to RLCT.
\begin{proof} \ \\
We have already proved $\lambda(K,\varphi) = \lambda(H',\varphi)$, hence we prove $\lambda(H',\varphi) = \lambda(H,\varphi)$.\\
Proof is almost same to previous proof.\\
Define $f_n, g_n$ as below.
\begin{eqnarray*}
f_n := \sum_{k=1}^H a_k b_k^{n} -\sum_{k=1}^{H^*} a_k^* {b_k^*}^{n} ,\\
g_n :=  \sum_{k=1}^H a_k b_k^{n}e^{-b_k} -\sum_{k=1}^{H^*} a_k^* {b_k^*}^{n} e^{-b_k^*}.\\
\end{eqnarray*}
By using this notation, $H(w)$ and $H'(w)$ are simply written like $H(w) = \sum_{x \in [0:H+H^*-1]} f_n^2$, $H'(w) = \sum_{x \in [0:H+H^*-1]} g_n^2$. By Taylor expansion of $e^{-b_k}$, $g_n$ is written using $f_n$ as below.
\begin{eqnarray*}
g_n = \sum_{s=0}^\infty \frac {(-1)^s} {s! } f_{n+s}.
\end{eqnarray*}
Firstly, we prove $\lambda(H',\varphi) \le \lambda(H,\varphi)$.\\
$f_{n+s} =  \sum_{i=0}^{H+H^*-1} \mathcal{F}_i^{(n+s)}  f_{i} $
where the definition of $ \mathcal{F}_i^{(n)} $ is same to the $ \mathcal{F}_i^{(n)} $ definition in Lemma \ref{mixed:poisson:lm01} so $ |\mathcal{F}_i^{(n)} |<R^n$. Then
\begin{eqnarray*}
g_n& =& \sum_{s=0}^\infty \frac {(-1)^s} {s! } f_{n+s}\\
& = & \sum_{s=0}^\infty \frac {(-1)^s} {s! }\left ( \sum_{i=0}^{H+H^*-1} \mathcal{F}_i^{(n+s)}  f_{i} \right ) \\
&=&\sum_{i=0}^{H+H^*-1}  \left ( \sum_{s=0}^\infty \frac {(-1)^s} {s! }  \mathcal{F}_i^{(n+s)} \right )  f_{i} .
\end{eqnarray*}
All that is left is proving $\sum_{s=0}^\infty \frac {(-1)^s} {s! }  \mathcal{F}_i^{(n+s)}$ is a bounded function. From $ |\mathcal{F}_i^{(n)} |<R^n$,
\[ \left | \sum_{s=0}^\infty \frac {(-1)^s} {s! }  \mathcal{F}_i^{(n+s)} \right | \le \sum_{s=0}^\infty \frac {R^s} {s! }  R^{n}
\le   R^n e^{R}
\] 
and $R$ is a constant value so $\sum_{s=0}^\infty \frac {(-1)^s} {s! }  \mathcal{F}_i^{(n+s)}$ is a bounded function on $W$.
It holds for all $g_n$. By RLCT's bounded function property, $\lambda(H',\varphi) \le \lambda(H,\varphi)$ follows.
Next, we prove $\lambda(H,\varphi) \le \lambda(H',\varphi)$.\\
The function $f_n$ is written using $g_n$ like
\begin{eqnarray*}
f_n &=& \sum_{k=1}^H a_k b_k^{n} -\sum_{k=1}^{H^*} a_k^* {b_k^*}^{n}\\
&=& \sum_{k=1}^H a_k b_k^{n} e^{-b_k} \left ( \sum_{j=n}^\infty \frac {b_k^{j-n}} {(j-n)!} \right ) -\sum_{k=1}^{H^*} a_k^* {b_k^*}^{n} e^{-b_k^*}\left ( \sum_{j=n}^\infty \frac {{b_k^*}^{j-n}} {(j-n)!} \right )\\
&=&  \sum_{j=n}^\infty \frac {1} {(j-n)!} \left ( \sum_{k=1}^H a_k b_k^{j} e^{-b_k}  -\sum_{k=1}^{H^*} a_k^* {b_k^*}^{j} e^{-b_k^*}   \right )=\sum_{j=n}^\infty \frac {g_j} {(j-n)!} .
\end{eqnarray*}
Then, $g_j$ is written like $g_j = \sum_{i=0}^{H+H^*-1} \mathcal{F}_i^{(j)}  g_{i} $ in the same way as Lemma \ref{mixed:poisson:lm01}. We can prove $\lambda(H,\varphi) \le \lambda(H',\varphi)$ using the similar arguments of proving $\lambda(H',\varphi) \le \lambda(H,\varphi)$.
\end{proof}
\end{lm.}
\subsection{Multidimensional case}
The multidimensional case is proved in the same way as the 1 dimensional case.  In concrete,  use Lemma \ref{lm3} instead of Lemma \ref{lm2}.
\begin{lm.}
$K'(w) = \sum_{x} \left ( p(x|w)-q(x) \right )^2$ is equivalent to
\begin{eqnarray*}
H(w) = \sum_{{\bf x} \in [0:H+H^*-1]^M} \left ( \sum_{k=1}^H a_k  {\bf b}_k^{\bf x} -\sum_{k=1}^{H^*} a_k^* {{\bf b^*}_k}^{\bf x} \right )^2
\end{eqnarray*}
with respect to RLCT.
\end{lm.}
In 2008, Aoyagi named the singularity of $H(w)$ Vandermode matrix type singularity\cite{aoyagi2}.
If there exists the restriction such that $a_k ,a_k^* \ge 0, b_k, b_k^* >0, \sum_{i=1}^H a_k = 1, \sum_{i=1}^{H^*} a_k^* = 1 $, we call singularity of $H(w)$ Simplex Vandermode matrix type singularity.
\section{Affine Variety of Simplex Vandermode Matrix Type Singularity}
The ideal generated by simplex Vandermode matrix type is
\[
 I_{\bf Po} := \bigg \langle \left \{ \sum_{k=1}^H a_k {\bf b}_k^{\bf x} - \sum_{k=1}^{H^*} a_k^* {\bf b}_k^{* \bf x} \right \}_{{\bf x} \in [0:H^*+H-1]^M} \bigg \rangle.
 \]
In this section, we study the properties of affine variety of simplex Vandermode matrix type.\\
Let $V( I_{\bf Po} )$ donate $\{ w \in W \mid \forall f \in  I_{\bf Po}, \left ( f(w) = 0 \right ) \}$.
\subsection{Preparation}
Consider the polynomial set $G_H = \left \{
\sum_{k=1}^H c_k s_k^x
\right \}_{x \in \mathbb{Z}_{\ge 0 }} $ where $\{c_k\}_{k \in [1:H]}$ in a set of variables and $\{ s_k \}_{k \in [1:H]}$ is constant values such that
\begin{eqnarray}
\forall k \in [1:H],\ s_k \in \mathbb{R}_{>0}  \label{eq:simple_case:01},
\end{eqnarray}
\begin{eqnarray}
\forall k, k' \in [1:H], \left (k \neq k' \Rightarrow s_k \neq s_{k'} \right ) \label{eq:simple_case:02}.
\end{eqnarray}

\begin{lm.} \label{lm:simple_case:01}
\[
\left \{ {\bf c} \in \mathbb{R}^H \mid \forall f \in G_H , \ f({\bf c}) = 0 \right \} = \{ {\bf 0} \}.
\]
\end{lm.}
\begin{proof}
\ \\
We resubscribe a new index of $s_k$ in descending order of $s_{i}$. That is, $ s_{1} > s_{2} > \cdots > s_{H} > 0$.
We prove Lemma by induction.
Firstly, prove $c_1 = 0$.\\
In the case $H=1$, it is trivial (think the case $x=0$).
In the case $H>1$,
for all $x \in  Z_{\ge 0 } $,
\begin{eqnarray*}
\sum_{k=1}^H c_k s_k^x = 0 \iff  c_1 =  \sum_{k=2}^H c_k \left ( \frac {s_k} {s_1} \right )^x .
\end{eqnarray*}
Since $0 \le \frac {s_k} {s_1} < 1$, $c_1 = \lim_{x \to \infty } \sum_{k=2}^H c_k \left ( \frac {s_k} {s_1} \right )^x = 0$.\\
Substitute $c_1 = 0$ and repeat the same operation to $c_2, c_3, \cdots, c_{H-1}$. We can prove $c_k = 0$ for all $k$ except $H$.
Then, $c_H = 0$ follows from $\sum_{k=1}^H c_k = 0$.
\end{proof}
In a multidimensional case,
the polynomial set $G_H = \left \{
\sum_{k=1}^H c_k {\bf s}_k^{\bf x}
\right \}_{{\bf x} \in \mathbb{Z}_{\ge 0 }^M} $ where $\{c_k\}_{k \in [1:H]}$ are variables, ${\bf x}$ is a multi-index and ${\bf s}_k$ is a M dimensional numerical vector such that 
\begin{eqnarray}
\forall k \in [1:H],\ {\bf s}_k \in \mathbb{R}^M \label{eq:simple_case:03},
\end{eqnarray}
\begin{eqnarray}
\forall k, k' \in [1:H], \left (k \neq k' \Rightarrow {\bf s}_k \neq {\bf s}_{k'} \right ) \label{eq:simple_case:04}.
\end{eqnarray}
Let $ s_{k,i}$ denote the i-th element of ${\bf s}_k$.
\begin{lm.}
\[
\left \{ {\bf c} \in \mathbb{R}^H \mid \forall f \in G_H , \ f({\bf c}) = 0 \right \} = \{ {\bf 0} \}.
\]\label{LinearIndepend:lm1}
\end{lm.}

\begin{proof}
\  \\
Define $f({\bf s'}_k, {\bf x'})$ such that $f({\bf s'}_k, {\bf x'}) := {\bf s'}_k ^ {\bf x'}$. Define a finite set $S_1$ as below.
\[
S_1 := \{ s_{k,1} \mid k \in [1:H] \}.
\]
Next, define $Inv(s_{k,1})$ given $s_{k,1}$ as below.
\[
Inv(s_{k,1}) := \{ k' \in [1:H] \mid s_{k',1} = s_{k,1}\}.
\]
Then, $\sum_{k=1}^H c_k f({\bf s'}_k, {\bf x'}) s_{k,1}^{ x_1}$ satisfies
\begin{eqnarray*}
\sum_{k=1}^H c_k f({\bf s'}, {\bf x'}) s_{k,1}^{ x_1} = \sum_{s_{k,1} \in S_1} \left ( \sum_{k \in Inv(s_{k,1})} c_k f({\bf s'}_k, {\bf x'})  \right )s_{k,1}^{ x_1}.
\end{eqnarray*}
Using the similar approach of Lemma \ref{lm:simple_case:01}, for all $s_{k,1} \in S_1$, we can prove $\sum_{k \in Inv(s_{k,1})} c_k f({\bf s'}_k, {\bf x'})   = 0$. Repeat the same operation to $\sum_{k \in Inv(s_{k,1})} c_k f({\bf s'}_k, {\bf x'})   = 0$ and then
\[
\sum_{k \in Inv({\bf s}_k)} c_k    = 0,
\]
where $Inv({\bf s}_k) := \cap_{i \in [1,M]} Inv(s_{k, i})$.\\
By condition (\ref{eq:simple_case:04}), $Inv({\bf s}_k) =\{ {\bf s}_k\}$. Thus, $c_k = 0$ for all $k \in [1:H]$.
\end{proof}
\subsection{Affine Variety}
We focus the structure of $V( I_{\bf Po} )$ using the method introduced in the previous subsection.
We use the notation as below.\\
\begin{eqnarray*}
Inv_{k,m} &:= &\{ k' \in [1:H] \mid b_{k', m} = b_{k, m}^* \},\\
Inv_k &:= &\cap_{m \in [1:M]} Inv_{k,m},\\
Inv_0& :=& [1:H] \setminus \left ( \cup_{k=1}^{H^*} Inv_k \right ).
\end{eqnarray*}
In this part, we state the properties of $V( I_{\bf Po} )$ which are easily proved.
\begin{itemize}
\item The affine variety $V( I_{\bf Po} )$ is equal to $V \left ( \bigg \langle \left \{ \sum_{k=1}^H a_k {\bf b}_k^{\bf x} - \sum_{k=1}^{H^*} a_k^* {\bf b}_k^{* \bf x} \right \}_{{\bf x} \in [0: H^*+H-1]^M} \bigg \rangle \right )$.
\item  The polynomial $\sum_{k=1}^H a_k {\bf b}_k^{\bf x} - \sum_{k=1}^{H^*} a_k^* {\bf b}_k^{* \bf x}$ is transformed to $   \sum_{k \in Inv_0} a_k {\bf b}_k^{\bf x}   + \sum_{i =1}^{H^* } \left (- a_i^* +  \sum_{k \in Inv_i} a_k \right ) {\bf b}_i^{* \bf x} $. By Lemma \ref{LinearIndepend:lm1} and  Condition $a_i \ge 0$, $({\bf a,b}) \in V( I_{\bf Po} ) \Rightarrow - a_i^* +  \sum_{j \in Inv_i} a_j  = 0$ for all $i\in [1:H^*]$ and $a_j = 0$ for all $j\in Inv_0$.
\item If $({\bf a,b}) \in V( I_{\bf Po} )$, then $Inv_i \neq \phi $ for all $i$ (if $Inv_i = \phi$, $- a_i^* +  \sum_{j \in Inv_i} a_j  = - a_i^* = 0$ contradicts the condition $a_i^*>0$.).
\end{itemize}
To summarize, the necessary and sufficient condition of $({\bf a,b}) \in V( I_{\bf Po} )$ is  $Inv_i \neq \phi $ and $- a_i^* +  \sum_{j \in Inv_i} a_j  = 0$ for all $i\in [1:H^*]$ and $a_j = 0$ for all $j\in Inv_0$.

That is, if we want to examine of the structure of $V( I_{\bf Po} )$ or to derive RLCT on $V( I_{\bf Po} )$, it is efficient to divide $V( I_{\bf Po} )$ into cases with respect to the number of  each $Inv_i$.\\

\section{Learning Coefficient (RLCT) of A Poisson Mixture}
We restrict a parameter region to a sufficient small neighborhood of $V( I_{\bf Po} )$.  This is because the $n$ order of $Z_n$ on $K(w)>\epsilon$ is $o_p({\rm exp}(- \sqrt {n}))$ (it is much smaller than  that of $Z_n$ on $K(w) \le \epsilon$) , we can ignore the points which are outside of the neighborhood of $V(I_{\bf Po})$.\\

In this section, we prove
\begin{th.}
\label{maintheorem}
Learning coefficient(RLCT) of Simplex Vandermode matrix type singularity is
\[
\lambda (K(w), \varphi(w))  =  \begin{cases}
    \frac {3r + H -2} {4} & (M=1) \\
     \frac {M r + H -1} {2}  & (M>1)
  \end{cases}.
\]
where M is the dimension of data $x$, H is the number of components of $p(x|w)$ and  $r$ is the number of components of $q(x)$.
\end{th.}
We have used the symbol $H^*$ as the number of components of true model, however,  we use the the symbol $r$ instead of $H^*$ to improve visual intelligibility from this section.

In 2018, we derived that RLCT of a Poisson mixture which is equal to RLCT of simplex Vandermode matrix type singularity satisfies
\[
\lambda (K(w), \varphi(w))  \le   \begin{cases}
    \frac {3r + H -2} {4} & (M=1) \\
     \frac {M r + H -1} {2}  & (M>1)
  \end{cases}.
\]
\cite{sato1}. If we prove that the lower bound matches the upper bound, Theorem \ref{maintheorem} follows.

\subsection{Aoyagi's decomposition theorem}
In 2010, Aoyagi proved that Vandermode matrix can be locally divided into simple shapes without changing RLCT\cite{aoyagi4}.
Let $A(w) \simeq_{w^*} B(w)$ denote $A(w)$ is sufficiently close to $B(w)$ on a neighborhood of the point $w^* \in W$. Variables
${\bf a}$ is defined by
\begin{eqnarray*}
&{\bf a}& := \left ( \begin{array}{cccccccccccc}
      a_1^{(1)} & a_2^{(1)} & \cdots &a_{H_1}^{(1)} & a_{1}^{(2)} &\cdots &a_{H_{r'}}^{(r')} & - a_1^* & - a_2^* & \cdots &  - a_r^*
    \end{array}\right ).\\
\end{eqnarray*}
Assume that ${\bf b_1^*}, {\bf b_2^*}, \cdots, {\bf b_r^* }, {\bf C}^{(r+1)}, \cdots, {\bf C}^{(r')}$ are different real vector.\\
For ${\bf l} = (l_1, l_2, \cdots, l_M) \in [0:H+r-1 ]^M $, ${\bf b}_i^{\bf l}$ is defined by ${\bf b}_i^{\bf l} := \prod_{m=1}^M b_{im}^{l_m}$. Then ${\bf B}_{\bf l}, {\bf B}$ are defined by
\begin{eqnarray*}
{\bf B}_{\bf l} := \left(
    \begin{array}{cccccc}
	{\bf b}_1^{(1) \bf l}  \\
	{\bf b}_2^{(1) \bf l}  \\
	\vdots \\
	{\bf b}_{H_1}^{(1) \bf l}  \\
	{\bf b}_1^{(2) \bf l} \\
	\vdots  \\
	{\bf b}_{H_{r'}}^{(r') \bf l}  \\
	{\bf b}_1^{* \bf l}  \\
	\vdots  \\
	{\bf b}_r^{*\bf l}   \\
    \end{array}
  \right),\ \ 
&{\bf B}& := \left(
    \begin{array}{cccccc}
	{\bf B}_{\bf l}  
    \end{array}
  \right)_{{\bf l} \in  [0:H+r -1]^M}.
\end{eqnarray*}
Also ${\bf a}^{(j)}, {\bf B}^{(j)}$ are defined by
\begin{eqnarray*}
&{\bf a}^{(j)}& := \left ( \begin{array}{cccccccccccc}
      a_1^{(j)} & a_2^{(j)} & \cdots &a_{H_j}^{(j)} &   - a_j^*
    \end{array}\right ) \ \mbox{\rm  if  } j \le r ,\\
\end{eqnarray*}
\begin{eqnarray*}
&{\bf a}^{(j)}& := \left ( \begin{array}{cccccccccccc}
      a_1^{(j)} & a_2^{(j)} & \cdots &a_{H_j}^{(j)}
    \end{array}\right ) \ \mbox{\rm  if  } j > r, \\
\end{eqnarray*}
and
\begin{eqnarray*}
{\bf B}_{\bf l}^{(j)} := \left(
    \begin{array}{cccccc}
	{\bf b}_1^{(j) \bf l}  \\
	{\bf b}_2^{(j) \bf l}  \\
	\vdots \\
	{\bf b}_{H_j}^{(j) \bf l}  \\
	{\bf b}_j^{*\bf l}   \\
    \end{array}
  \right),\ \ 
&{\bf B}^{(j)}& := \left(
    \begin{array}{cccccc}
	{\bf B}_{{\bf l}}^{(j)}  
    \end{array}
  \right)_{{\bf l} \in [0:H_j ]^M } \mbox{\rm  if  } j \le r,
\end{eqnarray*}
\begin{eqnarray*}
{\bf B}_{\bf l}^{(j)} := \left(
    \begin{array}{cccccc}
	{\bf b}_1^{(j) \bf l}  \\
	{\bf b}_2^{(j) \bf l}  \\
	\vdots \\
	{\bf b}_{H_j}^{(j) \bf l}  \\
    \end{array}
  \right),\ \ 
&{\bf B}^{(j)}& := \left(
    \begin{array}{cccccc}
	{\bf B}_{\bf l}^{(j)}  
    \end{array}
  \right)_{{\bf l} \in [0:H_j-1 ]^M } \mbox{\rm  if  } j > r.
\end{eqnarray*}
\begin{th.}[Aoyagi, 2010]
We have 
\[
\norm{\bf a B}^2 =_{RLCT} \sum_{j=1}^{r'} \norm{ {\bf a}^{(j)} {\bf B }^{(j)}}^2 \mbox{  on  small neighborhood of } w^* ,
\]
where the neighborhood of $w^*$ satisfies
\begin{eqnarray*}
\sum_{i=1}^{H_1} a_i^{(1)} \simeq_{w^*} {a^*}_1,&\hspace{5mm}& {\bf b}_1^{(1)}, {\bf b}_2^{(1)},\cdots, {\bf b}_{H_1}^{(1)}\simeq_{w^*} {\bf {b^*}}_1,\\
\sum_{i=1}^{H_2} a_i^{(2)} \simeq_{w^*} {a^*}_2,&\hspace{5mm}& {\bf b}_1^{(2)}, {\bf b}_2^{(2)},\cdots,{\bf b}_{H_2}^{(2)}\simeq_{w^*} {\bf {b^*}}_2,\\
&\vdots&\\
\sum_{i=1}^{H_r} a_i^{(r)} \simeq_{w^*} {a^*}_r,&\hspace{5mm}& {\bf b}_1^{(r)}, {\bf b}_2^{(r)},\cdots, {\bf b}_{H_r}^{(r)} \simeq_{w^*} {\bf {b^*}}_r, \\
a_1^{(r+1)} ,a_2^{(r+1)}\cdots,a_{H_{r+1}}^{(r+1)}\simeq_{w^*} 0,&\hspace{5mm}&{\bf b}_1^{(r+1)}\simeq_{w^*}  {\bf C}^{(r+1)}, {\bf b}_2^{(r+1)}\simeq_{w^*}  {\bf C}^{(r+1)}, \cdots\\
&\hspace{5mm}&\cdots,  {\bf b}_{H_{r+1}}^{(r+1)} \simeq_{w^*}  {\bf C}^{(r+1)}, \hspace{5mm}  {\bf b}_{i}^{(r+1) }\not \simeq_{w^*} {\bf b}_j^* \mbox{ for all } i,j,\\
&\vdots&\\
a_1^{(r')} ,a_2^{(r')}\cdots,a_{H_{r'}}^{(r')}\simeq_{w^*} 0,&\hspace{5mm}&{\bf b}_1^{(r')}\simeq_{w^*} {\bf C}^{(r')},{\bf b}_2^{(r')}\simeq_{w^*}  {\bf C}^{(r')} ,\cdots\\
&\hspace{5mm}&\cdots , {\bf b}_{H_{r'}}^{(r')}\simeq_{w^*} {\bf C}^{(r')}, \hspace{5mm}  {\bf b}_{i}^{(r')} \not \simeq_{w^*} {\bf b}_j^* \mbox{ for all } i,j,
\end{eqnarray*}

\end{th.}

\subsection{Preparation}

Before deriving lower bounds of RLCT, we introduce useful properties of 3 type ideals on local region. \\
\ \\
{\bf Type 1}: Let a {\bf Type 1} ideal be $\langle
{\bf a}
{\bf B} \rangle_{w} $. Assume $b_1^*>0$, $\sum_{i=1}^H a_i$ is sufficiently close to $a^*$ and $b_i$ be  sufficiently close to $b_1^*$ for all $i\in [1:H]$.
\begin{eqnarray*}
&{\bf a}& := \left ( \begin{array}{cccccc}
      - a^* & a_1 & a_2 & \cdots & a_H 
    \end{array}\right ),\\
&{\bf B}& := \left(
    \begin{array}{cccccc}
      1 & {b_1^*} & \cdots &  {b_1^*}^{H} \\
	1 & {b_1} & \cdots & {b_1}^{H} \\
	\vdots & \vdots & \vdots & \vdots  \\
	1 & {b_H} & \cdots &  {b_H}^{H}
    \end{array}
  \right),
\end{eqnarray*}
where $\langle {\bf a}{\bf B} \rangle_{w}$ is the ideal $\langle {\bf a}{\bf B} \rangle$ whose parameter region is restricted to a neighborhood of $w$. In this part, $w$ can be chosen as any point which satisfies the assumption.\\
Transform ${\bf B}$ by the algorithm below.
\begin{algorithm}[H]
\begin{algorithmic}
\FOR {$j \leftarrow H+1$ to $2$ }
	\STATE ${\bf B}_{\cdot, j} \leftarrow {\bf B}_{\cdot, j} - b_1^* {\bf B}_{\cdot, j-1}$
\ENDFOR
\end{algorithmic}
\end{algorithm}  
This transformation is justified by RLCT's ideal invariance property.\\
\begin{eqnarray*}
&\ &\langle 
{\bf a}
{\bf B} \rangle_w \\
&=_{RLCT}&\left \langle 
{\bf a}
\left(
    \begin{array}{cccccc}
      1 & 0&0 & \cdots & 0\\
	1 &   b_1 - b_1^* &  b_1 (b_1 - b_1^*) &\cdots & {b_1}^{H-1} (b_1 - b_1^*) \\
	\vdots & \vdots & \vdots & \ddots & \vdots  \\
	1 & b_H - b_1^*& b_H (b_H - b_1^*)& \cdots & {b_H}^{H-1} (b_H - b_1^*)
    \end{array}
  \right) \right  \rangle_w .
\end{eqnarray*}
By  ${\bf B}_{\cdot, 3} \leftarrow {\bf B}_{\cdot, 3} - b_1^* {\bf B}_{\cdot, 2}$, it follows that
\begin{eqnarray*}
\langle 
{\bf a}
{\bf B} \rangle_w &=_{RLCT}&\left \langle 
{\bf a}
\left(
    \begin{array}{cccccc}
       1 & 0&0 & \cdots & 0\\
	 1& b_1 - b_1^* &  (b_1 - b_1^*)^2 &\cdots & {b_1}^{H-1} (b_1 - b_1^*) \\
	\vdots & \vdots & \vdots & \ddots & \vdots  \\
	1& b_H - b_1^*&  (b_H - b_1^*)^2& \cdots & {b_H}^{H-1} (b_H - b_1^*)
    \end{array}
  \right) \right  \rangle_w .
\end{eqnarray*}
Focusing on the third element,
\begin{eqnarray*}
\norm {
{\bf a}
\left(
    \begin{array}{c}
      0 \\
	 (b_1 - b_1^*)^2 \\
	 \vdots \\
	 (b_H - b_1^*)^2
    \end{array}
  \right)}^2
=  \left (\sum_{i=1}^H a_i (b_i - b_1^*)^2 \right )^2 .
\end{eqnarray*}
There exist positive constant values $A, B$ such that
\begin{eqnarray*}
A \sum_{i=1}^H  \left (a_i  (b_i - b_1^*)^2 \right )^2 \le \left (\sum_{i=1}^H a_i  (b_i - b_1^*)^2 \right )^2 \le B \sum_{i=1}^H  \left (a_i (b_i - b_1^*)^2 \right )^2.
\end{eqnarray*}
Therefore, we can substitute $  \{ a_i  (b_i - b_1^*)^2 \}_{i\in [1:H]}$ for $\sum_{i=1}^H a_i  (b_i - b_1^*)^2$($\because$RLCT's inequality: constant factor property).\\
Next, focusing on the second element,
\begin{eqnarray*}
\norm {
{\bf a}
\left(
    \begin{array}{c}
      0 \\
	  b_1 - b_1^* \\
	 \vdots \\
	b_H - b_1^*
    \end{array}
  \right)}^2
=  \left (\sum_{i=1}^H a_i (b_i - b_1^*) \right )^2 .
\end{eqnarray*}
Thus, 
\begin{eqnarray*}
&\ &\langle 
{\bf a}
{\bf B} \rangle_w \\
&=_{RLCT}&\left \langle 
\sum_{i=1}^H  a_i  - a^* ,\ \{ a_i  (b_i - b_1^*)^2 \}_{i\in [1:H]},\ \sum_{i=1}^H  a_i  (b_i - b_1^*)   \right  \rangle_w .
\end{eqnarray*}

{\bf Type 2}: Let {\bf Type 2} ideal be $\langle {\bf a}{\bf B} \rangle_w$. Assume $b_1^*>0$, $a_1$ is sufficiently close to $a^*$ and $b_1$ is sufficiently close to $b_1^*$.
\begin{eqnarray*}
{\bf a} := \left ( \begin{array}{cc}
      - a^* & a_1 
    \end{array}\right ), {\bf B} := \left(
    \begin{array}{cccccc}
      1 & {b_1^*}  \\
	1 & {b_1} 
    \end{array}
  \right).
\end{eqnarray*}
By ${\bf B}_{\cdot, 2} \leftarrow {\bf B}_{\cdot, 2} - b_1^* {\bf B}_{\cdot, 1}$, it follows that
\begin{eqnarray*}
&\ &\norm {
\left ( \begin{array}{cc}
      - a^* & a_1 
    \end{array}\right )
\left(
    \begin{array}{cccccc}
      1 &0  \\
	1 &  b_1 -b_1^*
    \end{array}
  \right)}^2 = \left ( - a^*  + a_1  \right )^2 + \left ( a_1(b_1 - b_1^*) \right )^2.
\end{eqnarray*}
Focusing on the last term, $a_1 $ is sufficiently close to $a_1^* >0$ and there exists $A,B>0$ such that
\begin{eqnarray*}
A (b_1 - b_1^*) \le a_1   (b_1 - b_1^*) \le  B  (b_1 - b_1^*).
\end{eqnarray*}
Thus, we can substitute $ (b_1 - b_1^*)^2$ for $\left ( a_1  (b_1 - b_1^*) \right )^2$ ($\because$ RLCT's inequality: constant factor property) .
As a result,  $\langle b_1 - b_1^*,\  a_1  - a^*  \rangle_w =_{RLCT} \langle {\bf a B} \rangle_w $.

{\bf Type 3}: Let a {\bf Type 3} ideal be $\langle {\bf a} {\bf B} \rangle_w $.\\
Assume $a_1,a_2, \cdots, a_H$ are sufficiently close to $0$ and  $b_1,b_2, \cdots, b_H >0$.
\begin{eqnarray*}
&{\bf a}& := \left ( \begin{array}{cccc}
      a_1 & a_2 & \cdots &a_H 
    \end{array}\right ),\\
&{\bf B}& := \left(
    \begin{array}{cccccc}
       1 & {b_1}&\cdots &  {b_1}^{H-1}  \\
	\vdots &\vdots & \ddots&\vdots \\
	1 & {b_H} & \cdots & {b_H}^{H-1}
    \end{array}
  \right),\\
\end{eqnarray*}
Focusing on the first element,
\begin{eqnarray*}
&\ &\norm {
{\bf a}
 \left(
    \begin{array}{cccccc}
      1 \\
	\vdots  \\
	1
    \end{array}
  \right)}^2 = \left ( \sum_{i=1}^H a_i \right ) ^2.
\end{eqnarray*}
 Since $ a_i  \ge0$,
\begin{eqnarray*}
 \sum_{i=1}^H  a_i ^2 \le \left ( \sum_{i=1}^H a_i  \right ) ^2 \le H^2 \sum_{i=1}^H  a_i^2.
\end{eqnarray*}
Therefore, $\langle \{ a_i \}_{i\in[1:H]} \rangle_w=_{RLCT} \langle {\bf a} {\bf B} \rangle_w$.

\subsection{Local RLCT in the case dimension 1}
Assume that
\begin{eqnarray*}
r' \ge r ,\ a_i^{(j)} \ge 0, b_i >0 \mbox{ for all } i\in [1:H_j],\ j \in [1:r'], 
\end{eqnarray*}
\begin{eqnarray*}
 a_i^* >0, b_i^* >0 \mbox{ for all } i\in [1:r],
\end{eqnarray*}
\begin{eqnarray*}
\sum_{j=1}^{r'} \sum_{i=1}^{H_j} a_i^{(j)} = 1,&\ & \sum_{i=1}^{r} a_i^* = 1.
\end{eqnarray*}
Also assume that $b_1^*,b_2^*,\cdots,b_r^*$ are distinct.
For $H_i \in \mathbb{Z}_{\ge 0}$, $H = \sum_{i=1}^{r'} H_i$ and $H_i \ge 1$ for all $i \in [1:r]$, we study a neighborhood defined by
\begin{eqnarray*}
\sum_{i=1}^{H_1} a_i^{(1)} \simeq_{w^*} {a^*}_1,&\hspace{5mm}& b_1^{(1)}, b_2^{(1)},\cdots, b_{H_1}^{(1)}\simeq_{w^*} {b^*}_1,\\
\sum_{i=1}^{H_2} a_i^{(2)} \simeq_{w^*} {a^*}_2,&\hspace{5mm}& b_1^{(2)}, b_2^{(2)},\cdots, b_{H_2}^{(2)}\simeq_{w^*} {b^*}_2,\\
&\vdots&\\
\sum_{i=1}^{H_r} a_i^{(r)} \simeq_{w^*} {a^*}_r,&\hspace{5mm}& b_1^{(r)}, b_2^{(r)},\cdots, b_{H_r}^{(r)}\simeq_{w^*} {b^*}_r ,\\
a_1^{(r+1)} ,a_2^{(r+1)}\cdots,a_{H_{r+1}}^{(r+1)}\simeq_{w^*} 0,&\hspace{5mm}& b_1^{(r+1)}\simeq_{w^*}  C^{(r+1)}, b_2^{(r+1)}\simeq_{w^*}  C^{(r+1)}, \cdots\\
&\hspace{5mm}&\cdots,  b_{H_{r+1}}^{(r+1)}\simeq_{w^*}  C^{(r+1)}, \hspace{5mm}  b_{i}^{(r+1)} \not \simeq_{w^*} b_j^* \mbox{ for all } i,j,\\
&\vdots&\\
a_1^{(r')} ,a_2^{(r')}\cdots,a_{H_{r'}}^{(r')}\simeq_{w^*} 0,&\hspace{5mm}& b_1^{(r')}\simeq_{w^*}  C^{(r')}, b_2^{(r')}\simeq_{w^*}  C^{(r')} ,\cdots\\
&\hspace{5mm}& \cdots ,  b_{H_{r'}}^{(r')}\simeq_{w^*}  C^{(r')}, \hspace{5mm}  b_{i}^{(r')} \not \simeq_{w^*} b_j^* \mbox{ for all } i,j,
\end{eqnarray*}
where $C^{(r+1)} , \cdots,  C^{(r')} $ are distinct constant values.\\
On the assumption, we derive RLCT of $||{\bf a}  {\bf B}||^2$ where ${\bf a} ,{\bf B}$ are defined as below. 
\begin{eqnarray*}
&{\bf a}& := \left ( \begin{array}{cccccccccccc}
      a_1^{(1)} & a_2^{(1)} & \cdots &a_{H_1}^{(1)} & a_{1}^{(2)} &\cdots &a_{H_{r'}}^{(r')} & - a_1^* & - a_2^* & \cdots &  - a_r^*
    \end{array}\right )\\
&{\bf B}& := \left(
    \begin{array}{cccccc}
	b_1^{(1)} & b_1^{(1)2} & \cdots & b_1^{(1) H+r} \\
	b_2^{(1)} & b_2^{(1)2} & \cdots & b_2^{(1) H+r} \\
	\vdots & \vdots &\ddots &\vdots \\
	b_{H_1}^{(1)} & b_{H_1}^{(1)2} & \cdots & b_{H_1}^{(1) H+r} \\
	b_1^{(2)} & b_1^{(2) 2} & \cdots & b_1^{(2) H+r} \\
	\vdots & \vdots &\ddots &\vdots \\
	b_{H_{r'}}^{(r')} & b_{H_{r'}}^{(r') 2} & \cdots & b_{H_{r'}}^{(r') H+r} \\
	b_1^* & b_1^{*2} & \cdots & b_1^{*H+r}  \\
	\vdots & \vdots &\ddots &\vdots \\
	b_r^* & b_r^{*2} & \cdots & b_r^{*H+r}  \\
    \end{array}
  \right).
\end{eqnarray*}
Firstly, Apply Aoyagi's decomposition theorem,
\[
\norm{\bf a B}^2 =_{RLCT} \sum_{j=1}^{r'} \norm{ {\bf a}^{(j)} {\bf B }^{(j)}}^2 \mbox{  on  small neighborhood of } w^* ,
\]
where ${\bf a}^{(j)}, {\bf B}^{(j)}$ are
\begin{eqnarray*}
&{\bf a}^{(j)}& := \left ( \begin{array}{cccccccccccc}
      a_1^{(j)} & a_2^{(j)} & \cdots &a_{H_j}^{(j)} &   - a_j^*
    \end{array}\right ) \mbox{ if } j \le r ,\\
\end{eqnarray*}
\begin{eqnarray*}
&{\bf a}^{(j)}& := \left ( \begin{array}{cccccccccccc}
      a_1^{(j)} & a_2^{(j)} & \cdots &a_{H_j}^{(j)}
    \end{array}\right ) \mbox{\rm  if  } j > r ,
\end{eqnarray*}
\begin{eqnarray*}
{\bf B}^{(j)} := \left(
    \begin{array}{cccccc}
1 & b_j^{*}  &\cdots & b_j^{* H_j}  \\
	1 & b_1^{(j) }&\cdots & b_1^{(j)H_j }\\
	1&b_2^{(j) } &\cdots  & b_2^{(j)H_j }\\
	\vdots &\vdots &\ddots &\vdots \\
	1& b_{H_j}^{(j) }&\cdots & b_{H_j}^{(j)H_j } 
    \end{array}
  \right),\ \ \mbox{\rm  if  } j \le r,
\end{eqnarray*}
\begin{eqnarray*}
{\bf B}^{(j)} := \left(
    \begin{array}{cccccc}
	1 & b_1^{(j) }&\cdots & b_1^{(j)H_j }\\
	1&b_2^{(j) } &\cdots  & b_2^{(j)H_j }\\
	\vdots & \vdots &\ddots &\vdots \\
	1& b_{H_j}^{(j) }&\cdots & b_{H_j}^{(j)H_j } 
    \end{array}
  \right),\ \ \mbox{\rm  if  } j > r.
\end{eqnarray*}

Next, divide into 3 cases by type of ${\bf B}^{(j)}$.\\
(I)  If $j \le r$ and $H_{j} > 1$,  
\begin{eqnarray*}
{\bf B}^{(j)}  := \left(
    \begin{array}{cccccc}
     1  &  b_{j}^{*} & \cdots &  b_{j}^{* H_{j}} \\
	 1 &  b_1^{(j)} & \cdots & b_1^{(j) H_{j}} \\
	\vdots & \vdots & \ddots & \vdots  \\
	1  &b_{H_{j}}^{(j) } & \cdots &b_{H_{j}}^{(j) H_{j}}
    \end{array}
  \right),
\end{eqnarray*}
This is {\bf Type 1} of the previous subsection, therefore
\begin{eqnarray*}
\langle {\bf a}^{(j)}
{\bf B}^{(j)} \rangle_w=_{RLCT}\left \langle 
\sum_{i=1}^H  a_i^{(j)}  - a^{*(j)} ,\ \{ a_i ^{(j)} (b_i^{(j)} - b_j^*)^2 \}_{i\in [1:H]},\ \sum_{i=1}^H  a_i^{(j)}  (b_i^{(j)} - b_j^*)   \right  \rangle_w.
\end{eqnarray*}
(II)  If $j \le r$ and $H_{j} = 1$,
  \begin{eqnarray*}
{\bf B}^{(j)}  := \left(
    \begin{array}{cccccc}
      1 &  b_{j}^* \\
	1 & b_1^{(j)}   
    \end{array}
  \right).
\end{eqnarray*}
This is {\bf Type 2} of the previous subsection, therefore
\begin{eqnarray*}
\langle {\bf a}^{(j)} {\bf B}^{(j)} \rangle =_{RLCT} \langle b_1^{(j)} - b_j^*,\  a_1^{(j)}  - a_j^*  \rangle.
\end{eqnarray*}
(III)  If $j > r$,  
\begin{eqnarray*}
{\bf B}^{(j)}  := \left(
    \begin{array}{cccccc}
	1& b_1^{(j) } & \cdots & b_1^{(j) H_{j}} \\
	\vdots & \vdots & \ddots & \vdots  \\
	1 & b_{H_{j}}^{(j) } & \cdots &b_{H_{j}}^{(j) H_{j}}
    \end{array}
  \right),
\end{eqnarray*}
This is {\bf Type 3} of the previous subsection, therefore
\begin{eqnarray*}
\langle {\bf a}^{(j) } {\bf B}^{(j)} \rangle=_{RLCT} \left \langle\{ a_i^{(j)} \}_{i\in [1:H_{j}]} \right \rangle.
\end{eqnarray*}

Apply this manipulation to all $j \in [1:r' ]$ and then re-subscript index $j \in [1:r]$ in descending order of $H_{j}$,
\begin{eqnarray*}
\norm {{\bf a}{\bf B}}^2&=_{RLCT}& \sum_{j=1}^{s} \left ( \left ( - a_{j}^*  + \sum_{i=1}^{H_{j}} a_i^{(j)} \right )^2  \right .\\
&\ &\left . + \sum_{i=1}^{H_{j}}  \left (a_i^{(j)}  (b_i^{(j)} - b_{j}^*)^2 \right )^2 + \left ( \sum_{i=1}^{H_{j}} a_i^{(j)}  (b_i^{(j)} - b_{j}^*) \right )^2 \right )\\
&\ & + \sum_{j=s+1}^{r}\left ( (b_1^{(j)} - b_{j}^*)^2+(a_1^{(j)} - a_{j}^* )^2 \right ) + \sum_{j = r+1}^{r'} \sum_{i=1}^{H_{j}} a_i^{(j)2},
\end{eqnarray*}
where $s$ is the size of the set $\{j \in \{1,2,\cdots, r\} \mid H_{j} \not = 1 \}$.\\

Lastly, we derive the lower bound of RLCT of $\norm {{\bf a}{\bf B}}^2$.\\
There exist $i(j) \in [1:H_{j} ]$ such that $a_{i(j)}^{(j)} > \frac {a_{j}^*} {H_{j}}$ for all $j \in [1:s ]$,
Convert $b_{i(j)}^{(j)}$ into $b_{i(j)}^{'(j)}$ by
\begin{eqnarray*}
b_{i(j)}^{'(j)} = \sum_{i=1}^{H_{j}}   a_i^{(j)}  (b_i^{(j)} - b_{j}^*) ,
\end{eqnarray*}
for all $j \in [ 1:s]$ (Jacobian determinant is $\left | \prod_{j=1}^{s} \frac {1} {a_{i(j)}^{(j)}} \right |$).\\
By $a_1^{(1)} :=1 - \sum_{i=2}^{H_{1}} a_i^{(1)} - \sum_{j=2}^{r'} \sum_{i=1}^{H_{j}} a_i^{(j)}$,
\begin{eqnarray*}
 \left ( - a_{1}^*  + a_1^{(1)} + \sum_{i=2}^{H_{1}} a_i^{(1)} \right )^2 &=&  \left ( \sum_{j=2}^r a_{j}^*  - \sum_{j=2}^{r'} \sum_{i=1}^{H_{j}} a_i^{(j)}\right )^2\\
&=&  \left ( \sum_{j=2}^r \left ( a_{j}^* - \sum_{i=1}^{H_{j}} a_i^{(j)} \right ) \right )^2\\
&\in &  \left \langle \left \{ - a_{j}^* + \sum_{i=1}^{H_{j}} a_i^{(j)} \right \}_{j \in [2:r]} \right \rangle.
\end{eqnarray*}
Thus, we can delete $\left ( - a_{1}^*  + a_1^{(1)} + \sum_{i=2}^{H_{1}} a_i^{(1)} \right )^2$ without changing RLCT.
Then, convert $a_{1}^{(j)}$ into $u^{(j)}$ by
\begin{eqnarray*}
u^{(j)} =- a_{j}^*  + \sum_{i=1}^{H_{j}} a_i^{(j)} ,
\end{eqnarray*}
for all $j \in [2: s ]$(Jacobian determinant is $\left | \prod_{j=2}^{s}\frac {1} { b_1^{(j)}} \right |$). 
Convert $a_{1}^{(j)}$ into $u^{(j)}$ by
\begin{eqnarray*}
u^{(j)} =a_1^{(j)}  - a_{j}^* ,
\end{eqnarray*}
for all $j \in [s+1: r]$ (Jacobian determinant is $1$). We obtain
\begin{eqnarray*}
\norm {{\bf a}{\bf B}}^2&=_{RLCT}&  \sum_{i=1}^{H_{1}}  \left (a_i^{(1)}  (b_i^{(1)} - b_{1}^*)^2 \right )^2 + \left (b_{i(1)}^{'(1)} \right )^2 \\
&\ &+ \sum_{j=2}^{s} \left (u^{(j)2}    + \sum_{i=1}^{H_{j}}  \left (a_i^{(j)}  (b_i^{(j)} - b_{j}^*)^2 \right )^2 + \left (b_{i(j)}^{'(j)} \right )^2 \right )\\
&\ & + \sum_{j=s+1}^{r}\left ( (b_1^{(j)} - b_{j}^*)^2+u^{(j)2} \right ) + \sum_{j = r+1}^{r'} \sum_{i=1}^{H_{j}} a_i^{(j)2}\\
&\ge_{RLCT}&  \sum_{i=1,i \not = i(1)}^{H_{1}}  \left (a_i^{(1)}  (b_i^{(1)} - b_{1}^*)^2 \right )^2 + \left (b_{i(1)}^{'(1)} \right )^2 \\
&\ &+ \sum_{j=2}^{s} \left ( u^{(j)2}   + \sum_{i=1, i \not = i(j)}^{H_{j}}  \left (a_i^{(j)}  (b_i^{(j)} - b_{j}^*)^2 \right )^2 + \left (b_{i(j)}^{'(j)} \right )^2 \right )\\
&\ & + \sum_{j=s+1}^{r}\left ( (b_1^{(j)} - b_{j}^*)^2+u^{(j)2} \right ) + \sum_{j = r+1}^{r'} \sum_{i=1}^{H_{j}} a_i^{(j)2}.
\end{eqnarray*}
RLCT of the last formula is
\begin{eqnarray*}
&\ &(H_1 - 1) \min \left ( \frac {1} {2}, \frac {1} {4} \right ) + \frac {1} {2} + \sum_{j=2}^s \left (\frac {1} {2} + \sum_{i=1, i \not = i(j)}^{H_{j}} \min \left ( \frac {1} {2}, \frac {1} {4} \right ) + \frac {1} {2} \right ) \\
&\ &+ \sum_{j=s+1}^r \left ( \frac {1} {2} + \frac {1} {2} \right ) + \sum_{j=r+1}^{r'} \sum_{i=1}^{H_{j}} \frac {1} {2}\\
&=& r - \frac {1} {2} +  \frac {\sum_{j=1}^r H_{j} - r} {4} + \frac {\sum_{j=r+1}^{r'} H_{j}} {2}.
\end{eqnarray*}
in conclusion, the lower bound of RLCT is $ r - \frac {1} {2} +  \frac {\sum_{j=1}^r H_{j} - r} {4} + \frac {\sum_{j=r+1}^{r'} H_{j}} {2}$.
\subsection{Local RLCT in the case dimension M}
Assume that
\begin{eqnarray*}
a_i^{(j)} \ge 0, {\bf b}_1^{(j)}, {\bf b}_2^{(j)}, \cdots, {\bf b}_{H_{j}}^{(j)} \in \mathbb{R}_{>0}^{M} \mbox{ for all } i\in [1:H],\ j \in [1: r'],
\end{eqnarray*}
\begin{eqnarray*}
a_i^* >0, {\bf b_1^*}, {\bf b_2^*}, \cdots, {\bf b_r^* } \in \mathbb{R}_{>0}^{M} \mbox{ for all } i\in [1: r],
\end{eqnarray*}
\begin{eqnarray*}
\sum_{j=1}^{r'} \sum_{i=1}^{H_j} a_i^{(j)} = 1,&\ & \sum_{i=1}^{r} a_i^* = 1,
\end{eqnarray*}
and ${\bf b_1^*}, {\bf b_2^*}, \cdots, {\bf b_r^* }, {\bf C}^{(r+1)}, \cdots, {\bf C}^{(r')}$ are different real vector.
For $H_i \in \mathbb{Z}_{\ge 0}$ and $H = \sum_{i=1}^{r'} H_i$, $H_i \ge 1$ for all $i \in [1:r]$, we study a neighborhood, 
\begin{eqnarray*}
\sum_{i=1}^{H_1} a_i^{(1)} \simeq_{w^*} {a^*}_1,&\hspace{5mm}& {\bf b}_1^{(1)}, {\bf b_2}^{(1)},\cdots, {\bf b}_{H_1}^{(1)}\simeq_{w^*} {\bf {b^*}}_1,\\
\sum_{i=1}^{H_2} a_i^{(2)} \simeq_{w^*} {a^*}_2,&\hspace{5mm}& {\bf b}_1^{(2)}, {\bf b}_2^{(2)},\cdots,{\bf b}_{H_2}^{(2)}\simeq_{w^*} {\bf {b^*}}_2,\\
&\vdots&\\
\sum_{i=1}^{H_r} a_i^{(r)} \simeq_{w^*} {a^*}_r,&\hspace{5mm}& {\bf b}_1^{(r)}, {\bf b}_2^{(r)},\cdots, {\bf b}_{H_r}^{(r)} \simeq_{w^*} {\bf {b^*}}_r, \\
a_1^{(r+1)} ,a_2^{(r+1)}\cdots,a_{H_{r+1}}^{(r+1)}\simeq_{w^*} 0,&\hspace{5mm}&{\bf b}_1^{(r+1)}\simeq_{w^*}  {\bf C}^{(r+1)}, {\bf b}_2^{(r+1)}\simeq_{w^*}  {\bf C}^{(r+1)}, \cdots\\
&\hspace{5mm}&\cdots,  {\bf b}_{H_{r+1}}^{(r+1)} \simeq_{w^*}  {\bf C}^{(r+1)}, \hspace{5mm}  {\bf b}_{i}^{(r+1) }\not \simeq_{w^*} {\bf b}_j^* \mbox{ for all } i,j,\\
&\vdots&\\
a_1^{(r')} ,a_2^{(r')}\cdots,a_{H_{r'}}^{(r')}\simeq_{w^*} 0,&\hspace{5mm}&{\bf b}_1^{(r')}\simeq_{w^*} {\bf C}^{(r')},{\bf b}_2^{(r')}\simeq_{w^*}  {\bf C}^{(r')} ,\cdots\\
&\hspace{5mm}&\cdots , {\bf b}_{H_{r'}}^{(r')}\simeq_{w^*} {\bf C}^{(r')}, \hspace{5mm}  {\bf b}_{i}^{(r')} \not \simeq_{w^*} {\bf b}_j^* \mbox{ for all } i,j.
\end{eqnarray*}
On the assumption, we calculate RLCT of $||{\bf a}  {\bf B}||^2$ where ${\bf a} ,{\bf B}$ are defined as below. 
\begin{eqnarray*}
&{\bf a}& := \left ( \begin{array}{cccccccccccc}
      a_1^{(1)} & a_2^{(1)} & \cdots &a_{H_1}^{(1)} & a_{1}^{(2)} &\cdots &a_{H_{r'}}^{(r')} & - a_1^* & - a_2^* & \cdots &  - a_r^*
    \end{array}\right ),\\
\end{eqnarray*}
${\bf l} = (l_1, l_2, \cdots, l_M) \in [0: H+r -1]^M $, $\ {\bf b}_i^{\bf l} := \prod_{m=1}^M b_{im}^{l_m}$,
\begin{eqnarray*}
{\bf B}_{\bf l} := \left(
    \begin{array}{cccccc}
	{\bf b}_1^{(1)\bf l}  \\
	{\bf b}_2^{(1)\bf l}  \\
	\vdots \\
	{\bf b}_{H_1}^{(1)\bf l}  \\
	{\bf b}_1^{(2)\bf l} \\
	\vdots  \\
	{\bf b}_{H_{r'}}^{(r')\bf l}  \\
	{\bf b}_1^{*\bf l}  \\
	\vdots  \\
	{\bf b}_r^{*\bf l}   \\
    \end{array}
  \right),\ \ 
&{\bf B}& := \left(
    \begin{array}{cccccc}
	{\bf B}_{\bf l}  
    \end{array}
  \right)_{{\bf l} \in [0: H+r -1]^M }.
\end{eqnarray*}
Then  it follows that
\begin{eqnarray*}
\norm{{\bf a B}}^2 &=& \sum_{{\bf l} \in [0:H+r-1]^M } \left ( \sum_{j=1}^{r} \left (\sum_{i=1}^{H_{j}} a_i^{(j)} {\bf b}_{i}^{(j) {\bf l}} - a_j^* {\bf b}_{j}^{\bf l} \right ) + \sum_{j=r+1}^{r'} \sum_{i=1}^{H_{j}} a_i^{(j)} {\bf b}_{i}^{(j) {\bf l}}  \right )^2.
\end{eqnarray*}
Firstly, Apply Aoyagi's decomposition theorem, 
\begin{eqnarray*}
\norm{{\bf a B}}^2 &=_{RLCT}& \sum_{j=1}^{r'} \norm{ {\bf a}^{(j)} {\bf B }^{(j)}}^2 \mbox{  on  small neighborhood of } w^* ,
\end{eqnarray*}
where ${\bf a}^{(j)}, {\bf B}^{(j)}$ are
\begin{eqnarray*}
&{\bf a}^{(j)}& := \left ( \begin{array}{cccccccccccc}
      a_1^{(j)} & a_2^{(j)} & \cdots &a_{H_j}^{(j)} &   - a_j^*
    \end{array}\right ) \mbox{ if } j \le r, \\
\end{eqnarray*}
\begin{eqnarray*}
&{\bf a}^{(j)}& := \left ( \begin{array}{cccccccccccc}
      a_1^{(j)} & a_2^{(j)} & \cdots &a_{H_j}^{(j)}
    \end{array}\right ) \mbox{\rm  if  } j > r, \\
\end{eqnarray*}
\begin{eqnarray*}
{\bf B}_{\bf l}^{(j)} := \left(
    \begin{array}{cccccc}
	{\bf b}_1^{(j) \bf l}  \\
	{\bf b}_2^{(j) \bf l}  \\
	\vdots \\
	{\bf b}_{H_j}^{(j) \bf l}  \\
	{\bf b}_j^{*\bf l}   \\
    \end{array}
  \right),\ \ 
&{\bf B}^{(j)}& := \left(
    \begin{array}{cccccc}
	{\bf B}_{{\bf l}}^{(j)}  
    \end{array}
  \right)_{{\bf l} \in [0:H_j ]^M } \mbox{\rm  if  } j \le r,
\end{eqnarray*}
\begin{eqnarray*}
{\bf B}_{\bf l}^{(j)} := \left(
    \begin{array}{cccccc}
	{\bf b}_1^{(j) \bf l}  \\
	{\bf b}_2^{(j) \bf l}  \\
	\vdots \\
	{\bf b}_{H_j}^{(j) \bf l}  \\
    \end{array}
  \right),\ \ 
&{\bf B}^{(j)}& := \left(
    \begin{array}{cccccc}
	{\bf B}_{\bf l}^{(j)}  
    \end{array}
  \right)_{{\bf l} \in [0:H_j-1 ]^M } \mbox{\rm  if  } j > r.
\end{eqnarray*}
Next, we simplify $\norm{{\bf a B}}^2$ as below.
\begin{eqnarray*}
\norm{{\bf a B}}^2 &=_{RLCT}& \sum_{j=1}^{r'} \norm{ {\bf a}^{(j)} {\bf B }^{(j)}}^2 \\
	&=&\sum_{j=1}^{r} \sum_{{\bf l} \in [0:H_j]^M } \left (  \sum_{i=1}^{H_{j}} a_i^{(j)} {\bf b}_{i}^{(j) {\bf l}} - a_j^* {\bf b}_{j}^{*\bf l}  \right )^2+\sum_{j=r+1}^{r'} \sum_{{\bf l} \in [0:H_j-1]^M } \left (  \sum_{i=1}^{H_{j}} a_i^{(j)} {\bf b}_{i}^{(j) {\bf l}}  \right )^2\\
&\ge_{RLCT}&\sum_{j=1}^{r}  \sum_{m=1}^M \sum_{k=0}^{H_j} \left (  \sum_{i=1}^{H_{j}} a_i^{(j)} {b}_{im}^{(j) k} - a_j^* { b}_{jm}^{*k}  \right )^2+\sum_{j=r+1}^{r'} \sum_{m=1}^M \sum_{k=0}^{H_j-1}\left (  \sum_{i=1}^{H_{j}} a_i^{(j)} { b}_{im}^{(j)k }  \right )^2.
\end{eqnarray*}
The polynomial $\sum_{k=0}^{H_j} \left (  \sum_{i=1}^{H_{j}} a_i^{(j)} {b}_{im}^{(j) k} - a_j^* { b}_{jm}^{*k}  \right )^2$ corresponds to {\bf Type 1} or {\bf Type 2} and $\sum_{k=0}^{H_j-1}\left (  \sum_{i=1}^{H_{j}} a_i^{(j)} {b}_{im}^{(j)k }  \right )^2$ corresponds to {\bf Type 3} in the previous subsection.\\
By the same way as 1 dimensional case,
\begin{eqnarray*}
\norm{{\bf a B}}^2 &\ge_{RLCT}&\sum_{j=1}^{s}  \sum_{m=1}^M \sum_{k=0}^{H_j} \left (  \sum_{i=1}^{H_{j}} a_i^{(j)} { b}_{im}^{(j) k} - a_j^* {b}_{jm}^{*k}  \right )^2+ \sum_{j=s+1}^{r}  \sum_{m=1}^M \sum_{k=0}^{H_j} \left (  a_1^{(j)} {b}_{1m}^{(j) k} - a_j^* {b}_{jm}^{*k}  \right )^2\\
&\ &+\sum_{j=r+1}^{r'} \sum_{m=1}^M \sum_{k=0}^{H_j-1}\left (  \sum_{i=1}^{H_{j}} a_i^{(j)} {\bf b}_{im}^{(j)k }  \right )^2\\
 &\ge_{RLCT}&\sum_{j=1}^{s}  \sum_{m=1}^M  \left ( \left ( -a_j^* + \sum_{i=1}^{H_j} a_i^{(j)} \right )^2+ \sum_{i=1}^{H_j} \left ( a_i^{(j)} \left ( b_{im}^{(j)} - b_{jm}^* \right )^2 \right )^2+ \left ( \sum_{i=1}^{H_{j}} a_i^{(j)} {b}_{im}^{(j) } - a_j^* {b}_{jm}^{*} \right )^2\right )\\
&\ &+ \sum_{j=s+1}^{r}  \sum_{m=1}^M  \left (  \left ( a_1^{(j)} - a_j^* \right )^2 +\left  ( {b}_{1m}^{(j) } -  {b}_{jm}^* \right )^2   \right )\\
&\ &+\sum_{j=r+1}^{r'} \sum_{i=1}^{H_{j}} a_i^{(j)2}\\
 &=_{RLCT}&\sum_{j=1}^{s}  \left ( -a_j^* + \sum_{i=1}^{H_j} a_i^{(j)} \right )^2 \\
&\ &  + \sum_{j=1}^{s}  \sum_{m=1}^M  \sum_{i=1}^{H_j} \left ( a_i^{(j)} \left ( b_{im}^{(j)} - b_{jm}^* \right )^2 \right )^2+  \sum_{j=1}^{s}  \sum_{m=1}^M\left ( \sum_{i=1}^{H_{j}} a_i^{(j)} {b}_{im}^{(j) } - a_j^* {b}_{jm}^{*} \right ) ^2 \\
&\ &+ \sum_{j=s+1}^{r}  \left ( a_1^{(j)} - a_j^* \right )^2  + \sum_{j=s+1}^{r}  \sum_{m=1}^M \left (  {b}_{1m}^{(j) } -  { b}_{jm}^* \right )^2  \\
&\ &+\sum_{j=r+1}^{r'} \sum_{i=1}^{H_{j}} a_i^{(j)2}.
\end{eqnarray*}

There exist $i(j) \in [1:H_{j} ]$ such that $a_{i(j)}^{(j)} > \frac {a_{j}^*} {H_{j}}$ for all $j \in [1: s ]$,
Convert $b_{i(j)m}^{(j)}$ into $b_{i(j)m}^{'(j)}$ by
\begin{eqnarray*}
b_{i(j)m}^{'(j)} = \sum_{i=1}^{H_{j}}    a_i^{(j)}  (b_{im}^{(j)} - b_{jm}^*) ,
\end{eqnarray*}
for all $j \in [1:s], m \in [1:M]$(Jacobian determinant is $\left |\prod_{m=1}^M \prod_{j=1}^{s} \frac {1} {a_{i(j)}^{(j)}} \right |$).\\
$a_1^{(1)} :=1 - \sum_{i=2}^{H_{1}} a_i^{(1)} - \sum_{j=2}^{r'} \sum_{i=1}^{H_{j}} a_i^{(j)}$,thus, we can delete $\left ( - a_{1}^*  +\sum_{i=1}^{H_{1}} a_i^{(1)} \right )^2$ without changing RLCT.\\\\
Convert $a_{1}^{(j)}$ into $u^{(j)}$ by
\begin{eqnarray*}
u^{(j)} =- a_{j}^* + \sum_{i=1}^{H_{j}} a_i^{(j)},
\end{eqnarray*}
for all $j \in [2:s]$(Jacobian determinant is $1$).
Then convert $a_{1}^{(j)}$ into $u^{(j)}$ by
\begin{eqnarray*}
u^{(j)} =a_1^{(j)}  - a_{j}^* ,
\end{eqnarray*}
for all $j \in [s+1:r]$(Jacobian determinant is $1$). It follows that
\begin{eqnarray*}
\norm{{\bf a T}}^2&\ge_{RLCT} &\sum_{j=2}^{s}  \left ( -a_j^* + \sum_{i=1}^{H_j} a_i^{(j)} \right )^2 \\
&\ &  + \sum_{j=1}^{s}  \sum_{m=1}^M  \sum_{i=1}^{H_j} \left ( a_i^{(j)} \left ( b_{im}^{(j)} - b_{jm}^* \right ) ^2 \right )^2+  \sum_{j=1}^{s}  \sum_{m=1}^M\left ( \sum_{i=1}^{H_{j}} a_i^{(j)} {\bf b}_{im}^{(j) } - a_j^* {\bf b}_{jm}^{*} \right ) ^2 \\
&\ &+ \sum_{j=s+1}^{r}  \left ( a_1^{(j)} - a_j^* \right )^2  + \sum_{j=s+1}^{r}  \sum_{m=1}^M \left (  {\bf b}_{1m}^{(j) } -  {\bf b}_{jm}^* \right )^2  \\
&\ &+\sum_{j=r+1}^{r'} \sum_{i=1}^{H_{j}} a_i^{(j)2}\\
&=_{RLCT}&  \sum_{m=1}^M \sum_{j=1}^s   b_{i(j)m}^{'(j)2}  + \sum_{j=2}^{r} u^{(j) 2} \\
&\ & +  \sum_{j=1}^{s} \sum_{\substack {i=1\\ i \not = i(j)} }^{H_{j}}a_i^{(j)2}  \left (  \sum_{m=1}^{M}  (b_{im}^{(j)} - b_{jm}^*)^4 \right )  + \sum_{j=s+1}^{r} \sum_{m=1}^{M} ({\bf b}_{1m}^{(j)} - {\bf b}_{jm}^*)^2  \\
&\ &+ \sum_{j=r+1}^{r'}  \sum_{i=1}^{H_{j}} a_i^{(j)2}.
\end{eqnarray*}
Finally, RLCT of the last formula is
\begin{eqnarray*}
&\ &\frac {Ms} {2} + \frac {r-1} {2} + \sum_{j=1}^s (H_{j} -1) \min \left ( \frac {1} {2}, \frac {M} {2} \right ) + \frac {M (r-s)} {2} + \frac {\sum_{j=r+1}^{r'} H_j} {2}\\
&=&\frac {Mr+H-1} {2} .
\end{eqnarray*}
\subsection{RLCT of Simplex Vandermode Matrix Type Singularity}
In the last 2 subsection, we derive local lower bound of RLCT of simplex Vandermode matrix type singularity.
Next, we derive RLCT with following the method discussed in Section 5.
In the case dimension = 1,
the local lower bound is $ r - \frac {1} {2} +  \frac {\sum_{j=1}^r H_{j} - r} {4} + \frac {\sum_{j=r+1}^{r'} H_{j}} {2}$ and the conditions of $\sum_{j=1}^r H_{j}$ and $\sum_{j=r+1}^{r'} H_{j}$ is $H_{j} > 1 \mbox{ for all } j \in [1:r]$ and
$\sum_{j=1}^{r'} H_{j} = H .$
Therefore, the lower bound of RLCT is  $ r - \frac {1} {2} +  \frac {H- r} {4}  = \frac {3r + H - 2} {4}$.
In the case dimension $> 1$,
The local lower bound of RLCT is  $ \frac {Mr+H-1} {2} $ regardless of the point on the affine variety.
As a result, the lower bounds of RLCT matches the upper bounds so we get the exact value of RLCT of a Poisson mixture.
\section{Conclusion}
We succeed in calculating RLCT of a simplex Vandermonde matrix type singularity which is equal to RLCT of a Poisson mixture.
This research clarifies that the mean of the Generalization error $G_n$ of a Poisson mixture is
\begin{eqnarray*}
\mathbb{E} [ G_n ] = L(w_0) + \frac {3r + H - 2} {4n}  + o \left ( \frac {1} {n} \right )  \ \ \mbox{dimension}=1,
\end{eqnarray*}
\begin{eqnarray*}
\mathbb{E} [ G_n ] = L(w_0) + \frac {Mr+H-1} {2n}  + o \left ( \frac {1} {n} \right )\ \ \mbox{dimension}>1,
\end{eqnarray*}
where $r = H^*$.

\section{Acknowledgments}
We would like to thank M.Aoyagi for useful discussions about her research in 2010\cite{aoyagi4}. Her theorem spectacularly simplifies our proof.

\end{document}